\documentclass[preprint,11pt]{elsarticle}

\usepackage{setspace}
\usepackage{url}
\usepackage{graphicx,graphics}
\usepackage{subfigure}
\usepackage{epsfig}
\usepackage{psfrag}
\usepackage{hyperref}
\usepackage{amsmath,amssymb,amsfonts,amsthm}
\usepackage{mathrsfs}
\usepackage{color}
\usepackage{float}
\usepackage{natbib}
\usepackage{tabularx}
\usepackage{multirow}

\theoremstyle{remark}
\newtheorem{thm}{Theorem}[section]
\newtheorem{rmk}[thm]{Remark}

\newcommand{\Eref}[1]{Equation (\ref{#1})}
\newcommand{\fref}[1]{Figure (\ref{#1})}

\newcommand{\frefs}[1]{Figures~\ref{#1}}

\newcommand{\KK}{\mathbf{K}}
\newcommand{\bfm}{\mathbf{M}}
\newcommand{\bn}{\mathbf{N}}
\newcommand{\DD}{\mathbf{D_b}}

\newcommand{\bveps}{\boldsymbol{\varepsilon}}

\newcommand{\BB}{\mathbf{B}}

\newcommand{\rmd}{\mathrm{d}}

\begin{document}

\begin{frontmatter}

\title{NURBS-based finite element analysis of functionally graded plates: static bending, vibration, buckling and flutter}

\author[weimar]{Navid Valizadeh \fnref{fn1}}
\author[india]{Sundararajan Natarajan \corref{cor1}\fnref{fn2}}
\author[cardiff]{Octavio A Gonzalez-Estrada\fnref{fn3}}
\author[weimar]{Timon Rabczuk \fnref{fn4}}
\author[seigen]{Tinh Quoc Bui \fnref{fn5}}
\author[cardiff]{St\'ephane PA Bordas \fnref{fn6}}

\cortext[cor1]{Corresponding author}

\address[weimar]{Institute of Structural Mechanics, Bauhaus-Universit\"{a}t
Weimar, Marienstra\ss{}e, Weimar}
\address[india]{School of Civil and Environmental Engineering, The University of New South Wales, Sydney, Australia}
\address[cardiff]{Institute of Mechanics and Advanced
Materials, Cardiff University, Cardiff, UK}
\address[seigen]{Chair of Structural Mechanics, Department of Civil Engineering, University of Siegen, Siegen, Germany}

\fntext[fn1]{\url navid.valizadeh@uni-weimar.de}
\fntext[fn2]{\url sundararajan.natarajan@gmail.com}
\fntext[fn3]{\url agestrada@gmail.com}
\fntext[fn4]{\url timon.rabczuk@uni-weimar.de}
\fntext[fn5]{\url bui-quoc@bauwesen.uni-siegen.de}
\fntext[fn6]{\url stephane.bordas@alum.northwestern.edu}

\begin{abstract}
In this paper, a non-uniform rational B-spline based iso-geometric finite element method is used to study the static and dynamic characteristics of functionally graded material (FGM) plates. The material properties are assumed to be graded only in the thickness direction and the effective properties are computed either using the rule of mixtures or by Mori-Tanaka homogenization scheme. The plate kinematics is based on the first order shear deformation plate theory (FSDT). The shear correction factors are evaluated employing the energy equivalence principle and a simple modification to the shear correction factor is presented to alleviate shear locking. Static bending, mechanical and thermal buckling, linear free flexural vibration and supersonic flutter analysis of FGM plates are numerically studied. The accuracy of the present formulation is validated against available three-dimensional solutions. A detailed numerical study is carried out to examine the influence of the gradient index, the plate aspect ratio and the plate thickness on the global response of functionally graded material plates.
	
\end{abstract}

\begin{keyword} 
	isogeometric analysis \sep functionally graded  \sep Reissner Mindlin plate \sep gradient index
	\sep Shear locking \sep finite elements \sep partition of unity \sep
	free vibration \sep buckling \sep flutter \sep boundary conditions	
\end{keyword}

\end{frontmatter}


\section{Introduction}
Since its introduction to decrease the thermal stresses in propulsion systems and in airframes for space application~\cite{koizumi1993}, functionally graded materials (FGMs) have led researchers to investigate the structural behaviour of such structures. FGMs are considered to be an alternative for certain class of aerospace structures exposed to high temperature environment. FGMs are characterized by a smooth transition from one material to another, thus circumventing high inter-laminar shear stresses and de-lamination that persists in laminated composites. Thus, for structural integrity, FGMs have advantages over the fiber-matrix composites.

\subsection{Background}
The investigation of the static and the dynamic behaviour of FGM structures is fairly well covered in the literature. Some of the important contributions are discussed here. Different plate theories, viz, FSDT~\cite{Reddy2000,Yang2002,Sundararajan2005}, second and other higher order accurate theory~\cite{Qian2004a,Ferreira2006,natarajanmanickam2012} have been used to describe plate kinematics. Existing approaches in the literature to study plate and shell structures made up of FGMs uses finite element method based on Lagrange basis functions~\cite{Reddy2000,ganapathiprakash2006,Sundararajan2005}, meshfree methods~\cite{Qian2004a,Ferreira2006}. All existing approaches show shear locking when applied to thin plates. Different techniques by which the locking phenomenon can be suppressed include:

\begin{itemize}
\item Retain the original interpolations and subsequently use an optimal integration rule for evaluating the bending and the shear terms;
\item Mixed interpolation technique~\cite{bathedvorkin1985};
\item Use field redistributed substitute shape functions~\cite{somashekarprathap1987,ganapathivaradan1991}; 
\item Discrete shear gap method~\cite{bletzingerbischoff2000};
\item Stabilized conforming nodal integration~\cite{wangchen2004}, i.e., strain smoothing, SFEM~\cite{nguyenrabczuk2008,nguyen-xuanrabczuk2008}.
\item Use $p$-adaptivity, for example Moving Least Square approximations~\cite{kanok-nukulchaibarry2001}.
\end{itemize}

He \textit{et al.,}~\citep{He2001} presented a finite element formulation based on thin plate theory for the vibration control of FGM plates with integrated piezoelectric sensors and actuators under mechanical load, whereas Liew \textit{et al.,}~\citep{Liew1994} have analyzed the active vibration control of plates subjected to a thermal gradient using shear deformation theory. The parametric resonance of FGM plates is discussed in ~\cite{Ng2000} by Ng \textit{et al.,} based on Hamilton's principle and the assumed mode technique. Yang and Shen~\citep{Yang2001,Yang2002} have analyzed the dynamic response of thin FGM plates subjected to impulsive loads using a Galerkin Procedure coupled with modal superposition methods, whereas, by neglecting the heat conduction effect, such plates and panels in thermal environments have been examined based on shear deformation with temperature dependent material properties~\cite{Yang2002}. The static deformation and vibration of FGM plates based on higher-order shear deformation theory is studied by Qian \textit{et al.,}~\citep{Qian2004a} using the meshless local Petrov-Galerkin method (MLPG) and Natarajan and Ganapathi~\cite{natarajanmanickam2012} using shear flexible elements. Matsunaga~\cite{matsunaga2008} presented analytical solutions for simply supported rectangular FGM plates based on second-order shear deformation theory, whereas, three dimensional solutions are proposed in~\citep{Vel2002,Vel2004} for vibrations of simply supported rectangular FGM plates. Reddy~\citep{Reddy2000} presented finite element solutions for the dynamic analysis of FGM plates and Ferreira \textit{et al.,}~\citep{Ferreira2006} performed dynamic analysis of FGM plates based on higher order shear and normal deformable plate theory using MLPG. Birman~\cite{birman1995} and Javaheri and Eslami~\cite{javaherieslami2002} have studied buckling of FGM plates subjected to in-plane compressive loading. Woo \textit{et al.,}~\cite{woomeguid2003} analyzed the thermo-mechanical postbuckling behaviour of plates and shallow cylindrical FGM panels using a classical theory. Ganapathi \textit{et al.,}~\cite{ganapathiprakash2006}, using a $\mathcal{C}^o$ shear flexible quadrilateral element, studied buckling of non-rectangular FGM plates under mechanical and thermal loads. Prakash and Ganapathi~\citep{Prakash2006} studied the linear flutter characteristics of FGM panels exposed to supersonic flow. Haddadpour \textit{et al.,}~\citep{Haddadpour2007} and Sohn and Kim~\citep{Sohn2008,Sohn2009} investigated the nonlinear aspects of flutter characteristics using the finite element method. FGM plates, like other plate structures, may develop flaws. Recently, Yang and Chen~\cite{Yang2010} and Kitipornchai \textit{et al.,}~\cite{Kitipornchai2009} studied the dynamic characteristics of FGM beams with an edge crack. Natarajan \textit{et al.,}~\cite{natarajanbaiz2011,natarajanbaiz2011a} and Baiz \textit{et al.,}~\cite{baiznatarajan2011} studied the influence of the crack length on the free flexural vibrations of FGM plates using the XFEM and smoothed XFEM, respectively.

\subsection{Approach}
The main objective of this paper is to investigate the potential of NURBS based iso-geometric finite element methods to study the static and dynamic characteristics of Reissner-Mindlin plates. The present formulation also suffers from shear locking when lower order NURBS functions are used as basis functions. da Vaiga \textit{et al.,}~\cite{veigabuffa2012} showed that the shear locking phenomena can be suppressed by using higher order NURBS functions. A similar approach was employed to suppress shear locking in the element-free Galerkin method~\cite{nukulchaibarry2001}. In this paper, we propose a simple technique to suppress shear locking, which relies on the introduction of an artificial shear correction factor ~\cite{kikuchiishii1999} when lower order NURBS basis functions are used. The drawback of this approach is that the shear correction factor is problem dependent.

\subsection{Outline}
The paper is organized as follows. A brief overview on functionally graded materials and Reissner-Mindlin plate theory is presented in the next section. Section \ref{bspline} presents an overview of NURBS basis functions and a simple correction to the shear terms to alleviate shear locking. The efficiency of the present formulation, numerical results and parametric studies are presented in Section \ref{numexample}, followed by concluding remarks in the last section.

\section{Theoretical Formulation} \label{theory}
\subsection{Functionally graded material}
A rectangular plate made of a mixture of ceramic and metal is considered with the coordinates $x,y$ along the in-plane directions and $z$ along the thickness direction (see \fref{fig:platefig}). The material on the top surface $(z=h/2)$ of the plate is ceramic rich and is graded to metal at the bottom surface of the plate $(z=-h/2)$ by a power law distribution. The effective properties of the FGM plate can be computed by using the rule of mixtures or by employing the Mori-Tanaka homogenization scheme.

Let $V_i (i=c,m)$ be the volume fraction of the phase material. The subscripts $c$ and $m$ refer to ceramic and metal phases, respectively. The volume fraction of ceramic and metal phases are related by $V_c + V_m = 1$ and $V_c$ is expressed as:

\begin{equation}
V_c(z) = \left( \frac{2z+h}{2h} \right)^n
\end{equation}

where $n$ is the volume fraction exponent $(n \geq 0)$, also known as the gradient index. The variation of the composition of ceramic and metal is linear for $n=$1, the value of $n=$ 0 represents a fully ceramic plate and any other value of $n$ yields a composite material with a smooth transition from ceramic to metal.

\subsubsection*{Rule of mixtures}
Based on the rule of mixtures, the effective property of a FGM is computed using the following expression:

\begin{equation}
P = P_c V_c + P_m V_m
\label{eqn:rulemix}
\end{equation}

\subsubsection*{Mori-Tanaka homogenization method}
Based on the Mori-Tanaka homogenization method, the effective Young's modulus and Poisson's ratio are computed from the effective bulk modulus $K$ and the effective shear modulus $G$ as~\cite{Sundararajan2005}

\begin{eqnarray}
\frac{K_{\rm eff}-K_m}{K_c-K_m} &=& \frac{V_c}{1 + V_m \frac{3(K_c-K_m)}{3K_m+4G_m}} \nonumber \\
\frac{G_{\rm eff}-G_m}{G_c-G_m} &=& \frac{V_c}{1 + V_m \frac{(G_c-G_m)}{(G_m+f_1)} }
\end{eqnarray}

where

\begin{equation}
f_1 = \frac{G_m (9K_m+8G_m)}{6(K_m+2G_m)}
\end{equation}

The effective Young's modulus $E_{\rm eff}$ and Poisson's ratio $\nu_{\rm eff}$ can be computed from the following relations:

\begin{eqnarray}
E_{\rm eff} = \frac{9 K_{\rm eff} G_{\rm eff}}{3K_{\rm eff} + G_{\rm eff}} \nonumber \\
\nu_{\rm eff} = \frac{3K_{\rm eff} - 2G_{\rm eff}}{2(3K_{\rm eff} + G_{\rm eff})}
\label{eqn:young}
\end{eqnarray}

The effective mass density $\rho$ is computed using the rule of mixtures. The effective heat conductivity $\kappa_{\rm eff}$ and the coefficient of thermal expansion $\alpha_{\rm eff}$ is given by:

\begin{eqnarray}
\frac{\kappa_{\rm eff} - \kappa_m}{\kappa_c - \kappa_m} = \frac{V_c}{1 + V_m \frac{(\kappa_c - \kappa_m)}{3\kappa_m}} \nonumber \\
\frac{\alpha_{\rm eff} - \alpha_m}{\alpha_c - \alpha_m} = \frac{ \left( \frac{1}{K_{\rm eff}} - \frac{1}{K_m} \right)}{\left(\frac{1}{K_c} - \frac{1}{K_m} \right)}
\label{eqn:thermalcondalpha}
\end{eqnarray}

\subsubsection*{Temperature distribution through the thickness}
The temperature variation is assumed to occur in the thickness direction only and the temperature field is considered to be constant in the $xy$-plane. In such a case, the temperature distribution along the thickness can be obtained by solving a steady state heat transfer problem:

\begin{equation}
-{d \over dz} \left[ \kappa(z) {dT \over dz} \right] = 0, \hspace{0.5cm} T = T_c ~\textup{at}~ z = h/2;~~ T = T_m ~\textup{at} ~z = -h/2
\label{eqn:heat}
\end{equation}

The solution of \Eref{eqn:heat} is obtained by means of a polynomial series~\cite{wu2004} as

\begin{equation}
T(z) = T_m + (T_c - T_m) \eta(z,h)
\label{eqn:tempsolu}
\end{equation}

where,

\begin{equation}
\begin{split}
\eta(z,h) = {1 \over C} \left[ \left( {2z + h \over 2h} \right) - {\kappa_{cm} \over (n+1)\kappa_m} \left({2z + h \over 2h} \right)^{n+1} + \right. \\ 
\left. {\kappa_{cm} ^2 \over (2n+1)\kappa_m ^2 } \left({2z + h \over 2h} \right)^{2n+1}
-{\kappa_{cm} ^3 \over (3n+1)\kappa_m ^3 } \left({2z + h \over 2h} \right)^{3n+1} \right. \\ + 
\left. {\kappa_{cm} ^4 \over (4n+1)\kappa_m^4 } \left({2z + h \over 2h} \right)^{4n+1} 
- {\kappa_{cm} ^5 \over (5n+1)\kappa_m ^5 } \left({2z + h \over 2h} \right)^{5n+1} \right] ;
\end{split}
\label{eqn:heatconducres}
\end{equation}

\begin{equation}
\begin{split}
C = 1 - {\kappa_{cm} \over (n+1)\kappa_m} + {\kappa_{cm} ^2 \over (2 n+1)\kappa_m ^2} 
- {\kappa_{cm} ^3 \over (3n+1)\kappa_m ^3} \\ + {\kappa_{cm} ^4 \over (4n+1)\kappa_m ^4}
- {\kappa_{cm} ^5\over (5n+1)\kappa_m ^5}
\end{split}
\end{equation}

\subsection{Reissner-Mindlin Plates}
\begin{figure}[htpb]
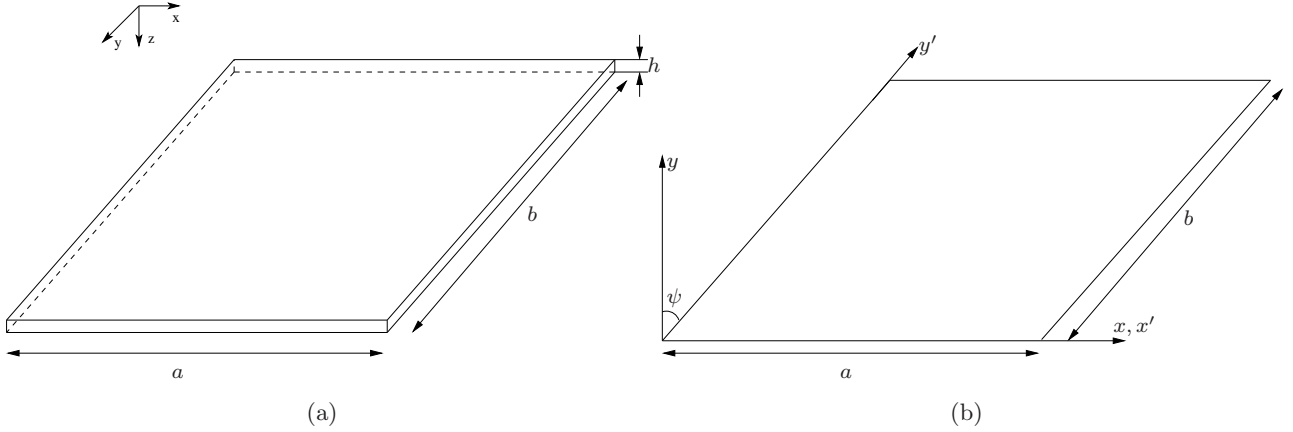

\centering
\subfigure[]{\input{plate.pstex_t}}
\subfigure[]{\input{skew.pstex_t}}
\caption{(a) coordinate system of a rectangular FGM plate, (b) Coordinate system of a skew plate}
\label{fig:platefig}
\end{figure}

The displacements $u,v,w$ at a point $(x,y,z)$ in the plate (see \fref{fig:platefig}) from the medium surface are expressed as functions of the mid-plane displacements $u_o,v_o,w_o$ and independent rotations $\theta_x,\theta_y$ of the normal in $yz$ and $xz$ planes, respectively, as:

\begin{eqnarray}
u(x,y,z,t) &=& u_o(x,y,t) + z \theta_x(x,y,t) \nonumber \\
v(x,y,z,t) &=& v_o(x,y,t) + z \theta_y(x,y,t) \nonumber \\
w(x,y,z,t) &=& w_o(x,y,t) 
\label{eqn:displacements}
\end{eqnarray}

where $t$ is the time. 

The strains in terms of mid-plane deformation can be written as:

\begin{equation}
\bveps  = \left\{ \begin{array}{c} \bveps_p \\ 0 \end{array} \right \}  + \left\{ \begin{array}{c} z \bveps_b \\ \bveps_s \end{array} \right\} 
\label{eqn:strain1}
\end{equation}

The midplane strains $\bveps_p$, bending strain $\bveps_b$, shear strain $\varepsilon_s$ in \Eref{eqn:strain1} are written as:

\begin{eqnarray}
\renewcommand{\arraystretch}{1.2}
\bveps_p = \left\{ \begin{array}{c} u_{o,x} \\ v_{o,y} \\ u_{o,y}+v_{o,x} \end{array} \right\}, \hspace{1cm}
\renewcommand{\arraystretch}{1.2}
\bveps_b = \left\{ \begin{array}{c} \theta_{x,x} \\ \theta_{y,y} \\ \theta_{x,y}+\theta_{y,x} \end{array} \right\}, \nonumber \\
\renewcommand{\arraystretch}{1.2}
\bveps_s = \left\{ \begin{array}{c} \theta _x + w_{o,x} \\ \theta _y + w_{o,y} \end{array} \right\}. \hspace{1cm}
\renewcommand{\arraystretch}{1.2}
\end{eqnarray}

where the subscript `comma' represents the partial derivative with respect to the spatial coordinate succeeding it. The membrane stress resultants $\bn$ and the bending stress resultants $\bfm$ can be related to the membrane strains, $\bveps_p$ and bending strains $\bveps_b$ through the following constitutive relations:

\begin{eqnarray}
\bn &=& \left\{ \begin{array}{c} N_{xx} \\ N_{yy} \\ N_{xy} \end{array} \right\} = \mathbf{A} \bveps_p + \BB \bveps_b \nonumber \\
\bfm &=& \left\{ \begin{array}{c} M_{xx} \\ M_{yy} \\ M_{xy} \end{array} \right\} = \BB \bveps_p + \DD \bveps_b 
\end{eqnarray}

where the matrices $\mathbf{A} = A_{ij}, \BB= B_{ij}$ and $\DD = D_{ij}; (i,j=1,2,6)$ are the extensional, bending-extensional coupling and bending stiffness coefficients and are defined as:

\begin{equation}
\left\{ A_{ij}, ~B_{ij}, ~ D_{ij} \right\} = \int_{-h/2}^{h/2} \overline{Q}_{ij} \left\{1,~z,~z^2 \right\}~dz
\end{equation}

Similarly, the transverse shear force $Q = \{Q_{xz},Q_{yz}\}$ is related to the transverse shear strains $\varepsilon_s$ through the following equation:

\begin{equation}
Q_{ij} = E_{ij} \varepsilon_s
\end{equation}

where $E_{ij} = \int_{-h/2}^{h/2} \overline{Q} \upsilon_i \upsilon_j~dz;~ (i,j=4,5)$ is the transverse shear stiffness coefficient, $\upsilon_i, \upsilon_j$ is the transverse shear coefficient for non-uniform shear strain distribution through the plate thickness. The stiffness coefficients $\overline{Q}_{ij}$ are defined as:

\begin{eqnarray}
\overline{Q}_{11} = \overline{Q}_{22} = {E(z) \over 1-\nu^2}; \hspace{1cm} \overline{Q}_{12} = {\nu E(z) \over 1-\nu^2}; \hspace{1cm} \overline{Q}_{16} = \overline{Q}_{26} = 0 \nonumber \\
\overline{Q}_{44} = \overline{Q}_{55} = \overline{Q}_{66} = {E(z) \over 2(1+\nu) }
\end{eqnarray}

where the modulus of elasticity $E(z)$ and Poisson's ratio $\nu$ are given by \Eref{eqn:young}. The thermal stress resultant $\bn^{\rm th}$ and the moment resultant $\bfm^{\rm th}$ are:

\begin{eqnarray}
\bn^{\rm th}&=& \left\{ \begin{array}{c} N^{\rm th}_{xx} \\ N^{\rm th}_{yy} \\ N^{\rm th}_{xy} \end{array} \right\} = \int\limits_{-h/2}^{h/2} \overline{Q}_{ij} \alpha(z,T) \left\{ \begin{array}{c} 1 \\ 1 \\ 0 \end{array} \right\} \Delta T(z)~ \rmd z \nonumber \\
\bfm^{\rm th} &=& \left\{ \begin{array}{c} M^{\rm th}_{xx} \\ M^{\rm th}_{yy} \\ M^{\rm th}_{xy} \end{array} \right\} = \int\limits_{-h/2}^{h/2} \overline{Q}_{ij} \alpha(z,T) \left\{ \begin{array}{c} 1 \\ 1 \\ 0 \end{array} \right\} \Delta T(z)~ z ~\rmd z \nonumber \\ 
\end{eqnarray}

where the thermal coefficient of expansion $\alpha(z,T)$ is given by \Eref{eqn:thermalcondalpha} and $\Delta T(z) = T(z)-T_o$ is the temperature rise from the reference temperature and $T_o$ is the temperature at which there are no thermal strains. The strain energy function $U$ is given by:

\begin{equation}
\begin{split}
U(\boldsymbol{\delta}) = {1 \over 2} \int_{\Omega} \left\{ \bveps_p^{\textup{T}} \mathbf{A} \bveps_p + \bveps_p^{\textup{T}} \mathbf{B} \bveps_b + 
\bveps_b^{\textup{T}} \mathbf{B} \bveps_p + \bveps_b^{\textup{T}} \mathbf{D} \bveps_b +  \bveps_s^{\textup{T}} \mathbf{E} \bveps_s \right\}~ \mathrm{d} \Omega
\end{split}
\label{eqn:potential}
\end{equation}

where $\boldsymbol{\delta} = \{u,v,w,\theta_x,\theta_y\}$ is the vector of the degree of freedom associated to the displacement field in a finite element discretization. Following the procedure given in~\cite{Rajasekaran1973}, the strain energy function $U$ given in~\Eref{eqn:potential} can be rewritten as:

\begin{equation}
U(\boldsymbol{\delta}) = {1 \over 2}  \boldsymbol{\delta}^{\textup{T}} \mathbf{K}  \boldsymbol{\delta}
\label{eqn:poten}
\end{equation}

where $\mathbf{K}$ is the linear stiffness matrix. The kinetic energy of the plate is given by:

\begin{equation}
T(\boldsymbol{\delta}) = {1 \over 2} \int_{\Omega} \left\{p (\dot{u}_o^2 + \dot{v}_o^2 + \dot{w}_o^2) + I(\dot{\theta}_x^2 + \dot{\theta}_y^2) \right\}~\mathrm{d} \Omega
\label{eqn:kinetic}
\end{equation}

where $p = \int_{-h/2}^{h/2} \rho(z)~dz, ~ I = \int_{-h/2}^{h/2} z^2 \rho(z)~dz$ and $\rho(z)$ is the mass density that varies through the thickness of the plate. When the plate is subjected to a temperature field, this in turn results in in-plane stress resultants, $\bn^{\rm th}$. The external work due to the in-plane stress resultants developed in the plate under a thermal load is given by:

\begin{equation}
\begin{split}
V(\boldsymbol{\delta}) = \int\limits_\Omega \left\{ \frac{1}{2} \left[ N_{xx}^{\rm th} w_{,x}^2 + N_{yy}^{\rm th} w_{,y}^2 + 2 N_{xy}^{\rm th}w_{,x}w_{,y}\right] + \right. \\ \left.
\frac{h^2}{24} \left[ N_{xx}^{\rm th} \left( \theta_{x,x}^2 + \theta_{y,x}^2 \right) + N_{yy}^2 \left( \theta_{x,y}^2 + \theta_{y,y}^2 \right) + 2 N_{xy}^{\rm th} \left( \theta_{x,x}\theta_{x,y} + \theta_{y,x}\theta_{y,y} \right) \right] \right\}~ d\Omega
\end{split}
\end{equation}

The governing equations of motion are obtained by writing the Lagrange equations of motion given by:

\begin{equation}
\frac{d}{dt} \left[ \frac{\partial(T-U)}{\partial \dot{\boldsymbol{\delta}}_i} \right] - \left[ \frac{\partial(T-U)}{\partial \boldsymbol{\delta}_i} \right] = 0, \hspace{0.5cm} i=1,2,\cdots,n
\end{equation}

The governing equations obtained using the minimization of total potential energy are solved using Galerkin finite element method. The finite element equations thus derived are:

{\bf Static bending}: 
\begin{equation}
\left( \KK + \KK_G \right) \boldsymbol{\delta} = \mathbf{F}
\label{eqn:staticgovern}
\end{equation}

{\bf Free vibration}: 
\begin{equation}
\bfm \ddot{\boldsymbol{\delta}} + \left( \KK + \KK_G \right) \boldsymbol{\delta} = \mathbf{0}
\label{eqn:freevibgovern}
\end{equation}

{\bf Buckling analysis}:

\paragraph*{Mechanical Buckling}
\begin{equation}
\left( \KK + \lambda \KK_G \right) \boldsymbol{\delta} = \mathbf{0}
\end{equation}

\paragraph*{Thermal Buckling}
\begin{equation}
\left( \KK + \Delta T \KK_G \right) \boldsymbol{\delta} = \mathbf{0}
\end{equation}

where $\boldsymbol{\delta}$ is the vector of degree of freedom associated to the displacement field in a finite element discretization, $\Delta T(=T_c - T_m)$ is the critical temperature difference, $\lambda$ is the critical buckling load and $\KK$, $\KK_G$ are the linear stiffness and geometric stiffness matrices, respectively. The critical temperature difference is computed using a standard eigenvalue algorithm. \\

{\bf Flutter analysis}: 
The work done by the applied non-conservative loads is:
\begin{equation}
W(\boldsymbol{\delta}) = \int_{\Omega} \Delta p w ~\rmd \Omega
\label{eqn:aerowork}
\end{equation}

where $\Delta p$ is the aerodynamic pressure. The aerodynamic pressure based on first-order, high Mach number approximation to linear potential flow is given by:

\begin{equation}
\Delta p = \frac{\rho_a U_a^2}{\sqrt{M_\infty^2 - 1}} \left[ \frac{\partial w}{\partial x} \cos \theta^\prime + \frac{\partial w}{\partial y} \sin \theta^\prime + \left( \frac{1}{U_a} \right) \frac{M_\infty^2 - 2}{M_\infty^2 - 1} \frac{\partial w}{\partial t} \right]
\label{eqn:aeropressure}
\end{equation}

where $\rho_a, U_a, M_\infty$ and $\theta^\prime$ are the free stream air density, velocity of air, Mach number and flow angle, respectively. The static aerodynamic approximation for Mach numbers between $\sqrt{2}$ and $2$ is~\citep{Birman1990}:

\begin{equation}
\Delta p = \frac{\rho_a U_a^2}{\sqrt{M_\infty^2 - 1}} \left[ \frac{\partial w}{\partial x} \cos \theta^\prime + \frac{\partial w}{\partial y} \sin \theta^\prime  \right]
\label{eqn:aeropressurestat}
\end{equation}

Substituting \Eref{eqn:poten} - (\ref{eqn:aerowork}) in Lagrange's equations of motion, the following governing equation is obtained:

\begin{equation}
\bfm \ddot{\boldsymbol{\delta}} + (\KK + \lambda \overline{\mathbf{A}}) \boldsymbol{\delta} = \mathbf{0}
\label{eqn:govereqn}
\end{equation}

After substituting the characteristic of the time function~\cite{Ganapathi1996} $\ddot{\boldsymbol{\delta}} = -\omega^2 \boldsymbol{\delta}$, the following algebraic equation is obtained:

\begin{equation}
\left[ \left( \KK + \lambda \overline{\mathbf{A}}\right) - \omega^2 \bfm\right] \boldsymbol{\delta} = \mathbf{0}
\label{eqn:finaldiscre}
\end{equation}

where $\KK$ is the stiffness matrix, $\bfm$ is the consistent mass matrix, $\lambda = \frac{\rho_a U_a^2}{\sqrt{M_\infty^2 - 1}}$, $\overline{\mathbf{A}}$ is the aerodynamic force matrix and $\omega$ is the natural frequency. When $\lambda = $0, the eigenvalue of $\omega$ is real and positive, since the stiffness matrix and mass matrix are symmetric and positive definite. However, the aerodynamic matrix $\overline{\mathbf{A}}$ is unsymmetric and hence complex eigenvalues $\omega$ are expected for $\lambda >$ 0. As $\lambda$ increases monotonically from zero, two of these eigenvalues will approach each other and become complex conjugates. In this study, $\lambda_{cr}$ is considered to be the value of $\lambda$ at which the first coalescence occurs.


\section{Non-Uniform Rational B-Splines}\label{bspline}

In this study, the finite element approximation uses NURBS basis function. We give here only a brief introduction to NURBS. More details on their use in FEM are given in~\cite{cottrellhughes,vinhphusimpson2012}. The key ingredients in the construction of NURBS basis functions are: the knot vector (a non decreasing sequence of parameter values, $\xi_i \le \xi_{i+1}, i = 0,1,\cdots,m-1$), the control points, $P_i$, the degree of the curve $p$ and the weight associated to a control point, $w$. The i$^{th}$ B-spline basis function of degree $p$, denoted by $N_{i,p}$ is defined as: 

\begin{eqnarray}
N_{i,0}(\xi) = \left\{ \begin{array}{cc} 1 & \textup{if} \hspace{0.2cm} \xi_i \le \xi \le \xi_{i+1} \\
0 & \textup{else} \end{array} \right. \nonumber \\
N_{i,p}(\xi) = \frac{ \xi- \xi_i}{\xi_{i+p} - \xi_i} N_{i,p-1}(\xi) + \frac{\xi_{i+p+1} - \xi}{\xi_{i+p+1}-\xi_{i+1}}N_{i+1,p-1}(\xi)
\end{eqnarray}

A $p^{th}$ degree NURBS curve is defined as follows:

\begin{equation}
\mathbf{C}(\xi) = \frac{\sum\limits_{i=0}^m N_{i,p}(\xi)w_i \mathbf{P}_i} {\sum\limits_{i=0}^m N_{i,p}(\xi)w_i}
\label{eqn:nurbsfunc1}
\end{equation}

\begin{figure}[htpb]
\centering
\includegraphics[scale=0.6]{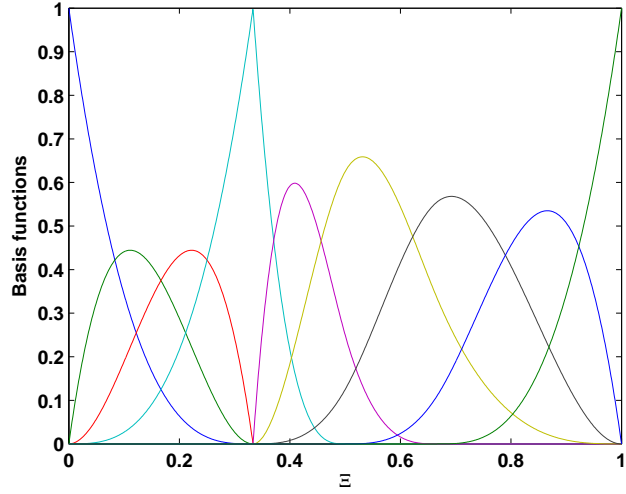}
\caption{non-uniform rational B-splines, order of the curve = 3}
\label{fig:nurbsplot}
\end{figure}

where $\mathbf{P}_i$ are the control points and $w_i$ are the associated weights. \fref{fig:nurbsplot} shows the third order non-uniform rational B-splines for a knot vector, $\Xi = \{0,~ 0,~ 0,~ 0,~ 1/3,~ 1/3,~ 1/3,~ 1/2,~ 2/3,~ 1,~ 1,~ 1,~ 1\}$. NURBS basis functions has the following properties: (i) non-negativity, (ii) partition of unity, $\sum\limits_i N_{i,p} = 1$; (iii) interpolatory at the end points. As the same function is also used to represent the geometry, the exact representation of the geometry is preserved. It should be noted that the continuity of the NURBS functions can be tailored to the needs of the problem. The B-spline surfaces are defined by the tensor product of basis functions in two parametric dimensions $\xi$ and $\eta$ with two knot vectors, one in each dimension as:

\begin{equation}
\mathbf{C}(\xi,\eta) = \sum_{i=1}^n\sum_{j=1}^m N_{i,p}(\xi)M_{j,q}(\eta) \mathbf{P}_{i,j}
\end{equation}

where $\mathbf{P}_{i,j}$ is the bidirectional control net and $N_{i,p}$ and $M_{j,q}$ are the B-spline basis functions defined on the knot vectors over an $m\times n$ net of control points $\mathbf{P}_{i,j}$. The NURBS surface is then defined by:

\begin{equation}
\mathbf{C}(\xi,\eta) = \frac{\sum_{i=1}^n\sum_{j=1}^m N_{i,p}(\xi)M_{j,q}(\eta) \mathbf{P}_{i,j} w_iw_j}{w(\xi,\eta)}
\label{eqn:nurbsfunc}
\end{equation}

where $w(\xi,\eta)$ is the weighting function. The displacement field within the control mesh is approximated by:

\begin{equation}
\{ u_o^e,v_o^e,w_o^e,\theta_x^e,\theta_y^e\} = \mathbf{C}(\xi,\eta) \{u_{oJ}, v_{oJ}, w_{oJ},\theta_{xJ},\theta_{yJ} \},
\end{equation} 

where $u_{oJ}, v_{oJ}, w_{oJ},\theta_{xJ},\theta_{yJ}$ are the nodal variables and $\mathbf{C}(\xi,\eta)$ are the basis functions given by \Eref{eqn:nurbsfunc}.




\section{Shear Locking} \label{shlock}
Transverse shear deformations are included in the formulation of Mindlin theory for thick plates. In Mindlin theory, the transverse normal to the mid surface of the plate before deformation remain straight but not necessarily normal to the mid surface after deformation. This relaxed the continuity requirement on the assumed displacement fields. But as the plate becomes very thin, care must be taken in not to violate the following relationship

\begin{equation}
\nabla w + \theta = 0
\end{equation}

i.e., the shear strain $\bveps_s$ must vanish in the domain as the thickness approaches zero. 

\subsection*{Artificial shear correction factor}
Lower order NURBS basis functions, like any other function, suffer from shear locking when applied to thin plates. Kikuchi and Ishii~\cite{kikuchiishii1999} introduced an artificial shear correction factor to suppress shear locking in 4-noded quadrilateral element. In this paper, we employ the same technique to suppress the shear locking syndrome in lower order NURBS basis functions. The modified shear correction factor is given by:

\begin{equation}
\upsilon_e = \upsilon \frac{ \left( \frac{h}{\beta l_e} \right)^2}{ \left( 1 + \left( \frac{h}{\beta l_e} \right)^{2n} \right)^{1/n}}
\end{equation}

where $n,\beta$ are positive integers, $l_e$ is the diameter of the element (=maximum or diagonal length of the element), $h$ is thickness of the element and $\upsilon$ is the shear correction factor. Here, the shear correction factor obtained based on energy equivalence principle as outlined in~\cite{singhaprakash2011,natarajanbaiz2011} is used.

\section{Numerical Examples} \label{numexample}
In this section, we present the static bending response, free vibration, buckling and flutter analysis of FGM plates using a NURBS based finite element method. The effect of various parameters, viz., material gradient index $n$, skewness of the plate $\psi$, the plate aspect ratio $a/b$, the plate thickness $a/h$ and boundary conditions on the global response is numerically studied. The top surface of the plate is ceramic rich and the bottom surface of the plate is metal rich. The material properties used for the FGM components are listed in Table \ref{tab:matprop}. 

\begin{table}[htbp]
\centering
\renewcommand{\arraystretch}{1.5}
\caption{Material properties.}
\begin{tabular}{lrrrr}
\hline
Property & Aluminum & Zirconia & Zirconia & Alumina \\
 & Al & ZrO$_2$-1 & ZrO$_2$-2 & Al$_2$O$_3$ \\
\hline
E (GPa) & 70 & 200 & 151 & 380 \\
$\nu$ & 0.3 & 0.3 & 0.3 & 0.3 \\
$\kappa$ W/mK & 204 & 2.09 & 2.09 & 10.4 \\
$\alpha$/$^\circ$C & 23 $\times$10$^{-6}$ & 10 $\times$10$^{-6}$ & 10 $\times$10$^{-6}$ & 7.2 $\times$10$^{-6}$ \\
$\rho$ kg/m$^3$ & 2707 & 5700 & 3000 & 3800 \\
\hline
\end{tabular}%
\label{tab:matprop}%
\end{table}%

\subsection*{Skew boundary transformation}
For skew plates supported on two adjacent edges, the edges of the boundary elements may not be parallel to the global axes $(x,y,z)$. In order to specify the boundary conditions on skew edges, it is necessary to use the edge displacements $(u_o^\prime,v_o^\prime,w_o^\prime)$, etc., in a local coordinate system $(x^\prime,y^\prime,z^\prime)$ (see \fref{fig:platefig}). The element matrices corresponding to the skew edges are transformed from global axes to local axes on which the boundary conditions can be conveniently specified. The relation between the global and local degrees of a particular node can be obtained through the following transformation~\cite{natarajanbaiz2011}

\begin{equation}
\boldsymbol{\delta} = \mathbf{L}_g \boldsymbol{\delta}^\prime 
\end{equation}

where $\boldsymbol{\delta}$ and $\boldsymbol{\delta}^\prime$ are the generalized displacement vector in the global and the local coordinate system, respectively. The nodal transformation matrix for a node $I$ on the skew boundary is given by

\begin{equation}
\mathbf{L}_g = \left[ \begin{array}{ccccc} \cos\psi & \sin\psi & 0 & 0 & 0 \\ -\sin\psi & \cos\psi & 0 & 0 & 0 \\ 0 & 0 & 1 & 0 & 0 \\ 0 & 0 & 0 & \cos\psi & \sin\psi \\ 0 & 0 & 0 & -\sin\psi & \cos\psi \end{array} \right]
\end{equation}

where $\psi$ defines the skewness of the plate.

\subsection{Static Bending}

Let us consider a Al/ZrO$_2$ FGM square plate with length-to-thickness $a/h=$ 5, subjected to a uniform load with fully simply supported (SSSS) and fully clamped (CCCC) boundary conditions. Four different values for the gradient index $(n=0,0.5,1,2)$ are considered in this study. The plate is modelled with 4, 8, 16, 24 and 32 control points per side. Tables \ref{tab:SSfgm1_static} and \ref{tab:CCfgm1_static} summarize the IGA results with quadratic, cubic and quartic NURBS elements for SSSS and CCCC boundary conditions. It can be seen that for all polynomial orders, the convergence of the results is quite fast. For cubic and quartic NURBS elements, the convergence is almost achieved with 16 control points per side. Table \ref{tab:ssss-cccc-valid} compares the results from the present formulation with other approaches available in the literature~\cite{gilhooleybatra2007,leezhao2009,nguyen-xuantran2012} and a very good agreement can be observed.


\begin{table}[htbp]
\centering
\renewcommand{\arraystretch}{1.5}
\caption{The normalized center deflection for fully simply supported Al/ZrO$_2$-1 FGM square plate with $a/h=$ 5, subjected to a uniformly distributed load $p$ using IGA. $^\ast$At least 5 control points per side are needed for 1 quartic NURBS element.}
\begin{tabular}{cccccc}
\hline
Method & Number of & \multicolumn{4}{c}{gradient index, $n$} \\
\cline{3-6}
& Control Points & 0 & 0.5 & 1 & 2 \\
\hline
\multirow{5}{*}{Quadratic} & 4 & 0.162098 & 0.218914 & 0.256018 & 0.293806 \\
& 8 & 0.171617 & 0.232392 & 0.271879 & 0.311459 \\
& 16 & 0.171649 & 0.232439 & 0.271935 & 0.311520 \\
& 24 & 0.171651 & 0.232441 & 0.271938 & 0.311523 \\
& 32 & 0.171651 & 0.232442 & 0.271938 & 0.311523 \\
\cline{2-6}
\multirow{5}{*}{Cubic} & 4 & 0.163329 & 0.220598 & 0.257991 & 0.296053 \\
& 8 & 0.171658 & 0.232452 & 0.271950 & 0.311536 \\
& 16 & 0.171651 & 0.232441 & 0.271938 & 0.311522 \\
& 24 & 0.171651 & 0.232442 & 0.271938 & 0.311523 \\
& 32 & 0.171651 & 0.232442 & 0.271938 & 0.311523 \\
\cline{2-6}
\multirow{5}{*}{Quartic} & 5$^\ast$ & 0.172910 & 0.234167 & 0.273959 & 0.313821 \\
& 8 & 0.171690 & 0.232495 & 0.272000 & 0.311594 \\
& 16 & 0.171651 & 0.232442 & 0.271938 & 0.311523 \\
& 24 & 0.171651 & 0.232442 & 0.271938 & 0.311523 \\
& 32 & 0.171651 & 0.232442 & 0.271938 & 0.311523 \\
\hline
\end{tabular}%
\label{tab:SSfgm1_static}%
\end{table}%

\begin{table}[htbp]
\centering
\renewcommand{\arraystretch}{1.5}
\caption{The normalized center deflection for fully clamped Al/ZrO$_2$-1 FGM square plate with $a/h=$ 5, subjected to a uniformly distributed load $p$ using IGA. $^\ast$At least 5 control points per side are needed for 1 quartic NURBS element.}
\begin{tabular}{cccccc}
\hline
Method & Number of & \multicolumn{4}{c}{gradient index, $n$} \\
\cline{3-6}
& Control Points & 0 & 0.5 & 1 & 2 \\
\hline
\multirow{5}{*}{Quadratic} & 4 & 0.052510 & 0.06818 & 0.079306 & 0.093421 \\
& 8 & 0.075831 & 0.101036 & 0.117946 & 0.136557 \\
& 16 & 0.076017 & 0.101305 & 0.118264 & 0.136905 \\
& 24 & 0.076024 & 0.101315 & 0.118276 & 0.136918 \\
& 32 & 0.076025 & 0.101316 & 0.118278 & 0.136920 \\
\cline{2-6}
\multirow{5}{*}{Cubic} & 4 & 0.056337 & 0.07342 & 0.085447 & 0.100403 \\
& 8 & 0.076031 & 0.101324 & 0.118287 & 0.136931 \\
& 16 & 0.076025 & 0.101316 & 0.118277 & 0.136921 \\
& 24 & 0.076025 & 0.101317 & 0.118278 & 0.136921 \\
& 32 & 0.076026 & 0.101317 & 0.118279 & 0.136921 \\
\cline{2-6}
\multirow{5}{*}{Quartic} & 5$^\ast$ & 0.077690 & 0.103620 & 0.120980 & 0.139973 \\
& 8 & 0.076092 & 0.101409 & 0.118387 & 0.137043 \\
& 16 & 0.076026 & 0.101317 & 0.118279 & 0.136921 \\
& 24 & 0.076026 & 0.101317 & 0.118279 & 0.136921 \\
& 32 & 0.076026 & 0.101317 & 0.118279 & 0.136921 \\
\hline
\end{tabular}%
\label{tab:CCfgm1_static}%
\end{table}%

\begin{table}[htbp]
\centering
\renewcommand{\arraystretch}{1.5}
\caption{The normalized center deflection for fully clamped Al/ZrO$_2$-1 FGM square plate with $a/h=$ 5, subjected to a uniformly distributed load $p$ using 13 $\times$13 control points for various boundary conditions. $^\ast$At least 5 control points per side are needed for 1 quartic NURBS element.}
\begin{tabular}{clcccc}
\hline
Method & Number of & \multicolumn{4}{c}{gradient index, $n$} \\
\cline{3-6}
& Control Points & 0 & 0.5 & 1 & 2 \\
\hline
\multirow{8}{*}{SSSS} & IGA-Quadratic & 0.1717 & 0.2324 & 0.2719 & 0.3115 \\
& IGA-Cubic & 0.1717 & 0.2324 & 0.2719 & 0.3115 \\
& IGA-Quartic & 0.1717 & 0.2324 & 0.2719 & 0.3115 \\
& NS-DSG3~\cite{nguyen-xuantran2012} & 0.1721 & 0.2326 & 0.2716 & 0.3107 \\
& ES-DSG3~\cite{nguyen-xuantran2012} & 0.1700 & 0.2296 & 0.2680 & 0.3066 \\
& MITC4~\cite{nguyen-xuantran2012} & 0.1715 & 0.2317 & 0.2704 & 0.3093 \\
& $kp-$Ritz~\cite{leezhao2009} & 0.1722 & 0.2403 & 0.2811 & 0.3221 \\
& MLPG~\cite{gilhooleybatra2007} & 0.1671 & 0.2505 & 0.2905 & 0.3280 \\
\cline{2-6}
\multirow{8}{*}{CCCC} & IGA-Quadratic & 0.0760 & 0.1013 & 0.1183 & 0.1369 \\
& IGA-Cubic & 0.0760 & 0.1014 & 0.1183 & 0.1369 \\
& IGA-Quartic & 0.0760 & 0.1014 & 0.1183 & 0.1369 \\
& NS-DSG3~\cite{nguyen-xuantran2012} & 0.0788 & 0.1051 & 0.1227 & 0.1420 \\
& ES-DSG3~\cite{nguyen-xuantran2012} & 0.0761 & 0.1013 & 0.1183 & 0.1370 \\
& MITC4~\cite{nguyen-xuantran2012} & 0.0758 & 0.1010 & 0.1179 & 0.1365 \\
& $kp-$Ritz~\cite{leezhao2009} & 0.0774 & 0.1034 & 0.1207 & 0.1404 \\
& MLPG~\cite{gilhooleybatra2007} & 0.0731 & 0.1073 & 0.1253 & 0.1444 \\
\hline
\end{tabular}%
\label{tab:ssss-cccc-valid}%
\end{table}%

\begin{figure}[htpb]
\centering
\subfigure[]{\includegraphics[scale=0.5]{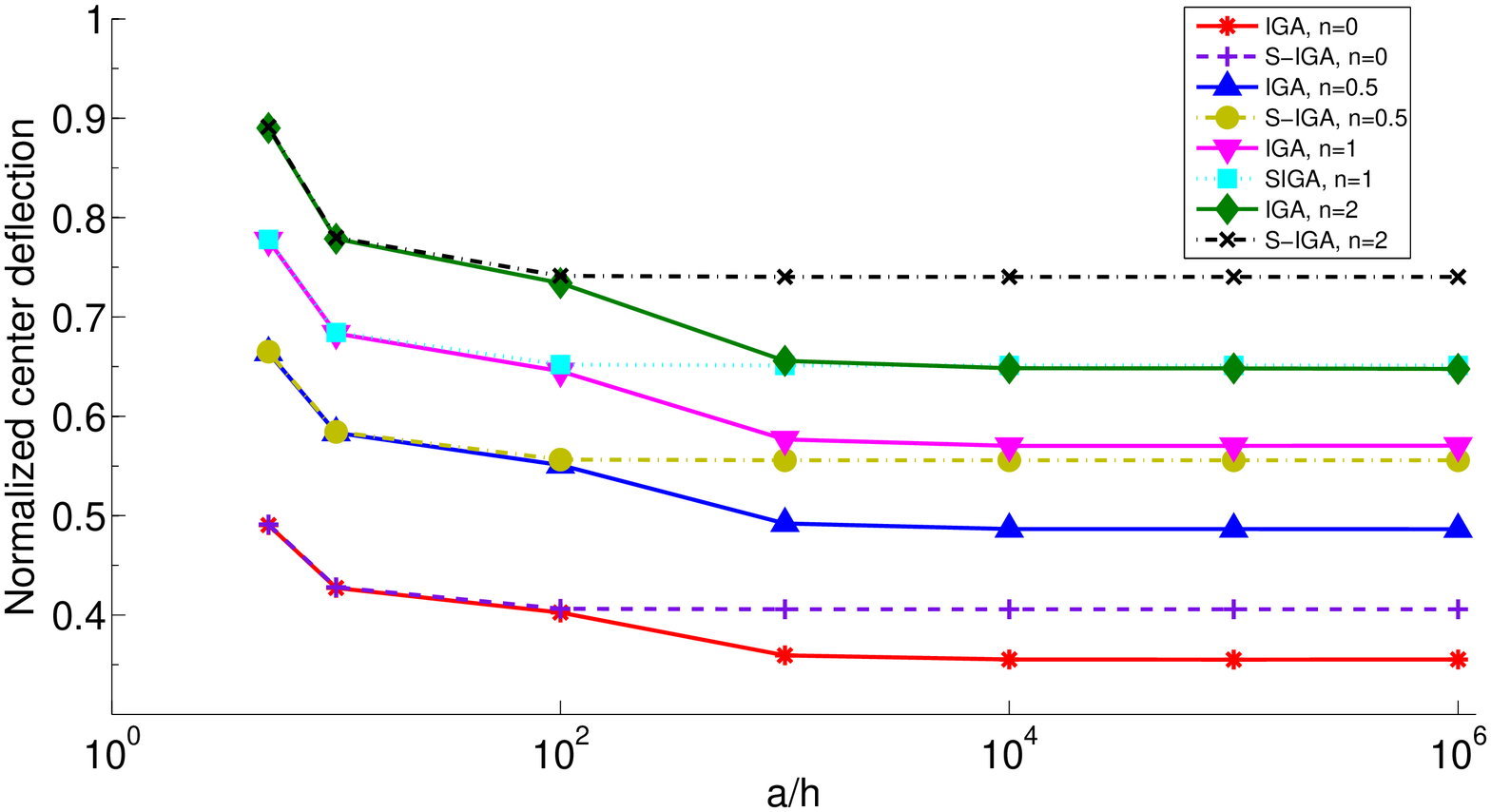}}
\subfigure[]{\includegraphics[scale=0.5]{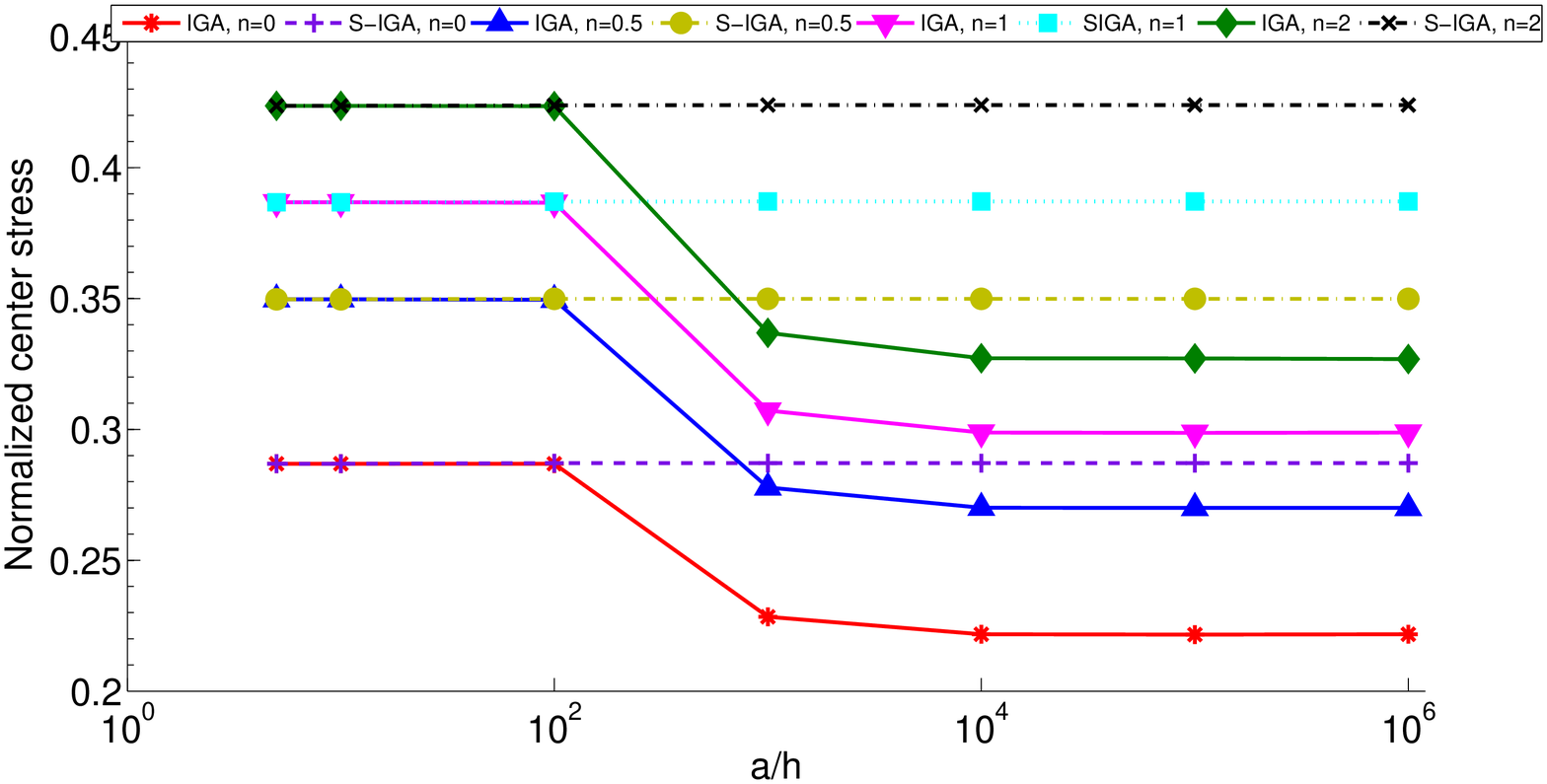}}
\caption{The normalized center deflection, $\overline{w}_c$ and the normalized axial stress $\overline{\sigma}_{xx}$ as a function of $a/h$ for a simply supported Al/Zr0$_2$-1 FGM square plate subjected to a uniform load, $p$: (a) normalized center deflection and (b) normalized center stress.}
\label{fig:SSS_defl_stres}
\end{figure}

\begin{figure}[htpb]
\centering
\subfigure[]{\includegraphics[scale=0.5]{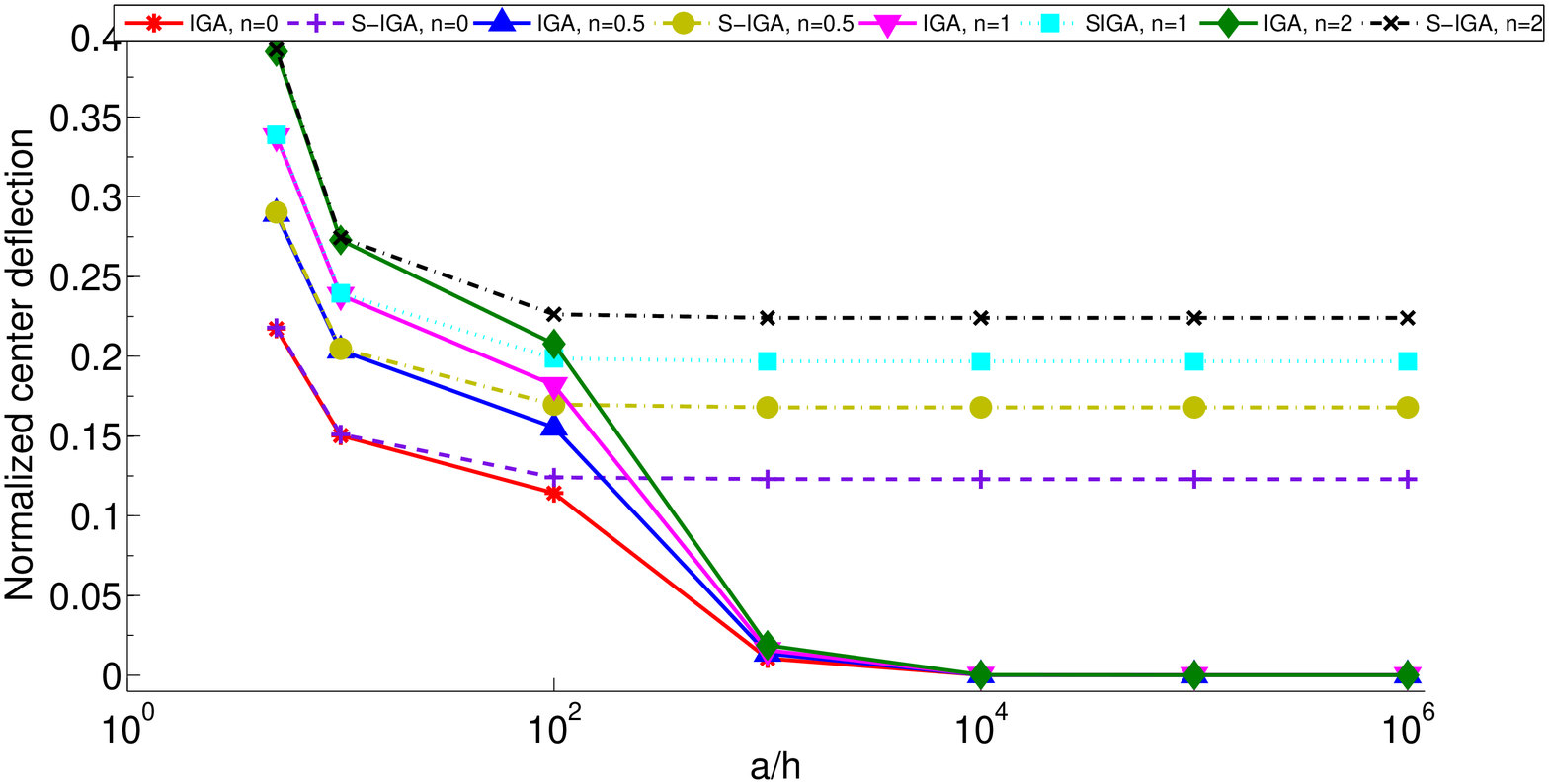}}
\subfigure[]{\includegraphics[scale=0.5]{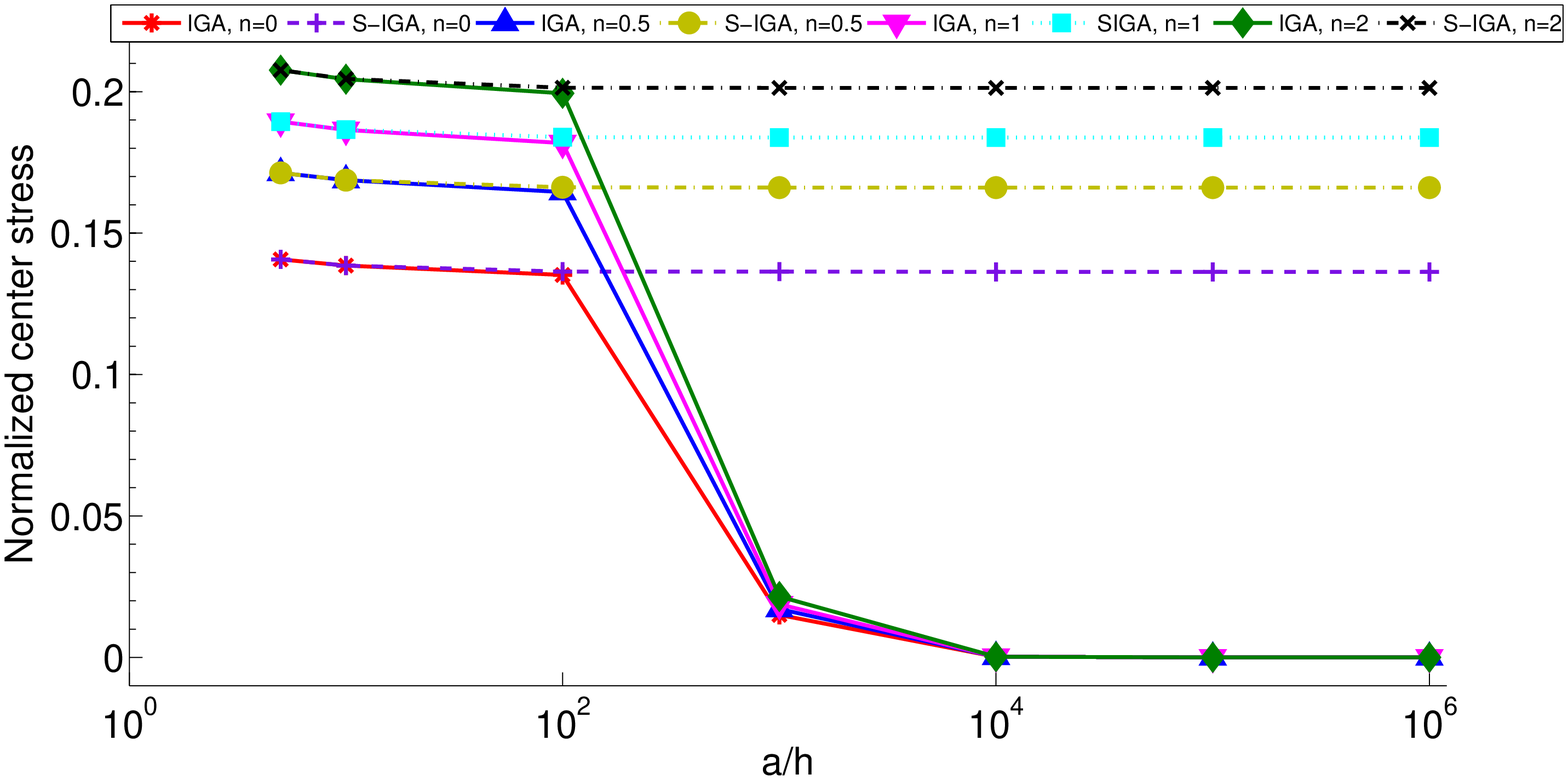}}
\caption{The normalized center deflection, $\overline{w}_c$ and the normalized axial stress $\overline{\sigma}_{xx}$ as a function of $a/h$ for a clamped Al/Zr0$_2$-1 FGM square plate subjected to a uniform load, $p$: (a) normalized center deflection and (b) normalized center stress.}
\label{fig:CCCC_defl_stres}
\end{figure}

Next, we illustrate the performance of the present isogeometric method for thin plate problems. A simply supported and a clamped Al/ZrO$_2$-1 square plates subjected to uniform load are considered, while the length-to-thickness ($a/h$) varies from 5 to 10$^6$ and the gradient index ranges from 0 to 2. Two individual approaches are employed: one applied the stabilization technique to eliminate shear locking named S-IGA and the other one, normal IGA, without considering any specific technique for shear locking. The plate is modelled using quadratic NURBS elements with 13$\times$13 control points. The normalized center deflection $\overline{w}_c = 100 w_c \frac{E_c h^3}{12(1-\nu^2)pa^4}$ and the normalized axial stress $\overline{\sigma}_{xx} = \sigma_{xx} \frac{h^2}{p a^2}$ at the top surface of the center of the plate for SSSS and CCCC boundary conditions are depicted in \frefs{fig:SSS_defl_stres} - \ref{fig:CCCC_defl_stres}, respectively. It is observed that IGA results are subjected to shear locking when the plate becomes thin ( $a/h>$ 100). However, the S-IGA results are almost independent of the length-to-thickness ratio for thin plates. The same observations also reported in~\cite{thainguyen-xuan2012} for laminated composite plates. The results using the S-IGA agree very well with those given in~\cite{nguyen-xuantran2012} using the NS-DSG3 element. 

\subsection*{FGM plates under thermo-mechanical loads}
\begin{figure}[htpb]
\centering
\includegraphics[scale=0.5]{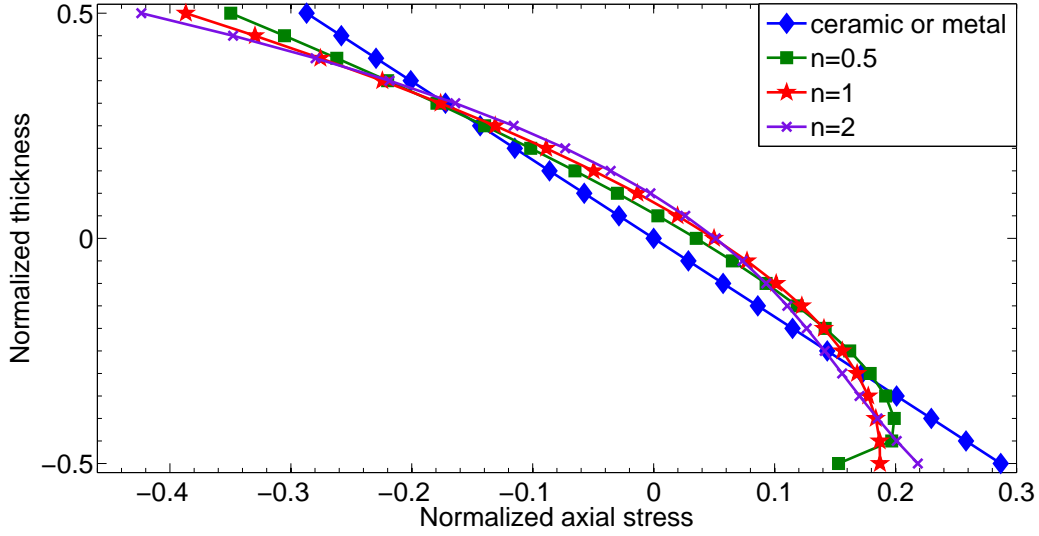}
\caption{The normalized axial stress $\overline{\sigma}_{xx} = \sigma_{xx} \frac{h^2}{p L^2}$ at points on the vertical line passing through the centroid of a simply supported square Al/Zr0$_2$ FGM plate subjected to uniform mechanical load $p$. }
\label{fig:SS_centerline_stres_fgm1}
\end{figure}

\begin{figure}[htpb]
\centering
\subfigure[]{\includegraphics[scale=0.5]{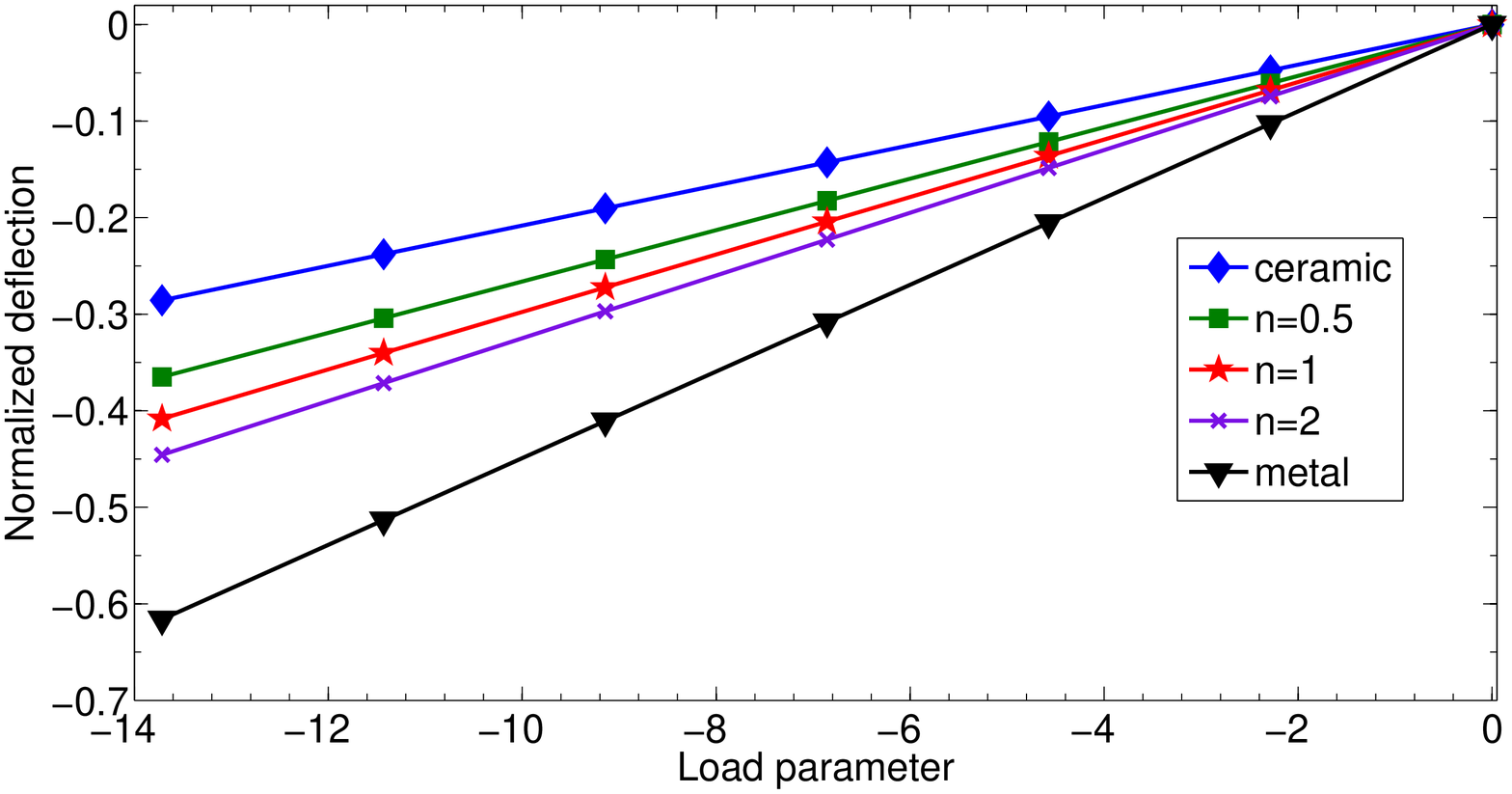}}
\subfigure[]{\includegraphics[scale=0.5]{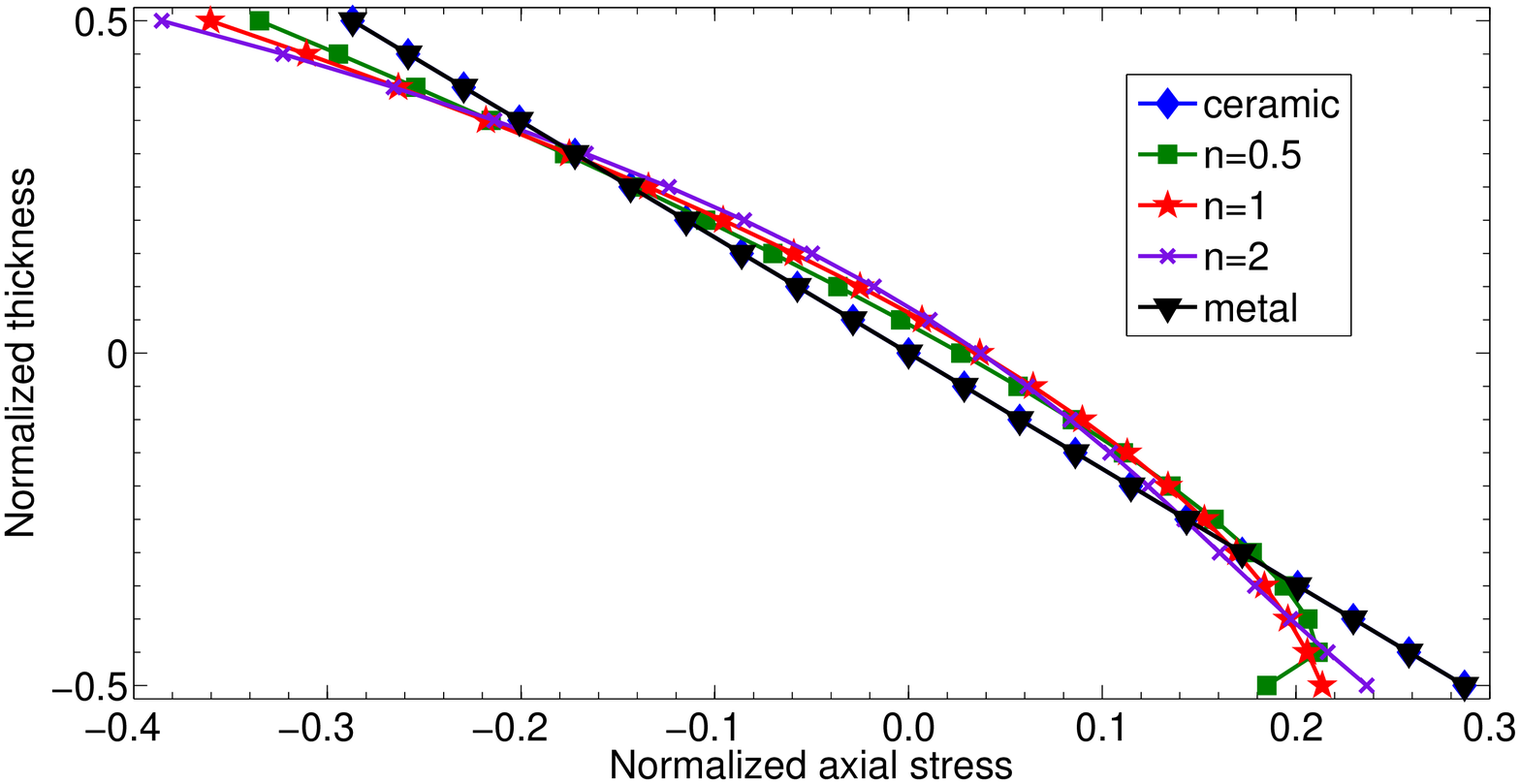}}
\caption{The normalized center deflection $\overline{w}_c$ and the normalized axial stress $\overline{\sigma}_{xx}$ of a simply supported square Al/Zr0$_2$ FGM plate subjected to uniform mechanical load $p$: (a)The normalized center deflection, $\overline{w}_c = w_c/h$with $\overline{p} = pa^4/(E_mh^4)$ and (b) The normalized axial stress $\overline{\sigma}_{xx} = \sigma_{xx} \frac{h^2}{p a^2}$ at points on the vertical line passing through the centroid with load $p$. }
\label{fig:SS_centerline_stres_fgm2}
\end{figure}

Here, the present isogeometric method is verified on FGM plates subjected to thermo-mechanical loads. A simply supported Al/ZrO$_2$ plate with aluminium at the bottom surface and zirconia at the top surface is considered. The plate with length $a=$ 0.2 m and thickness h=0.01 m is modelled employing quadratic NURBS elements. For comparison, a mesh with 13 control points per side is used. At first, we study the behaviour of SSSS Al/ZrO$_2$-I FGM plate under a uniform mechanical load. \fref{fig:SS_centerline_stres_fgm1} shows the distributions of the normalized axial stress through the thickness of the plate computed for different values of the gradient index $n$. The results are in excellent agreement with those given in~\cite{gilhooleybatra2007,nguyen-xuantran2012}. \fref{fig:SS_centerline_stres_fgm2} plots the central deflection $\overline{w}_c = 100 w_c \frac{E_c h^3}{12(1-\nu^2)pa^4}$ of the plate with respect to various load parameters, $\overline{p} = pa^4/(E_mh^4)$, given in the interval [-14, 0] for different values of the gradient index $n$. It can be seen that the central deflection of the plate linearly increases with respect to the load. It is also observed that the central deflection increases with the gradient index $k$. As expected, the metallic plate has the largest deflection, while the ceramic plate has the lowest. Note that the results match well with those given in~\cite{nguyen-xuantran2012,crocevenini2007}. \fref{fig:SS_centerline_stres_fgm2} shows the normalized axial stress at points on the vertical line passing through the centroid of simply supported Al/ZrO$_2$-2 FGM square plate subjected to a uniform mechanical load $p$. The results agree well with the results reported in~\cite{leezhao2009,nguyen-xuantran2012}. It is observed that for isotropic plates $(i.e, n = 0)$ the axial stress distribution is linear while it is nonlinear for FGM plates. For FGM plates, the magnitude of the axial stress at the bottom is less than the one at the top. The maximum compressive stress at the top surface of the plate has been obtained for the FGM plate with $n=$ 2, whilst the metallic only or ceramic only plates has the minimum tensile stress at the bottom surface.

\subsection*{FGM plates subjected to thermal loading}
\begin{figure}[htpb]
\centering
\includegraphics[scale=0.5]{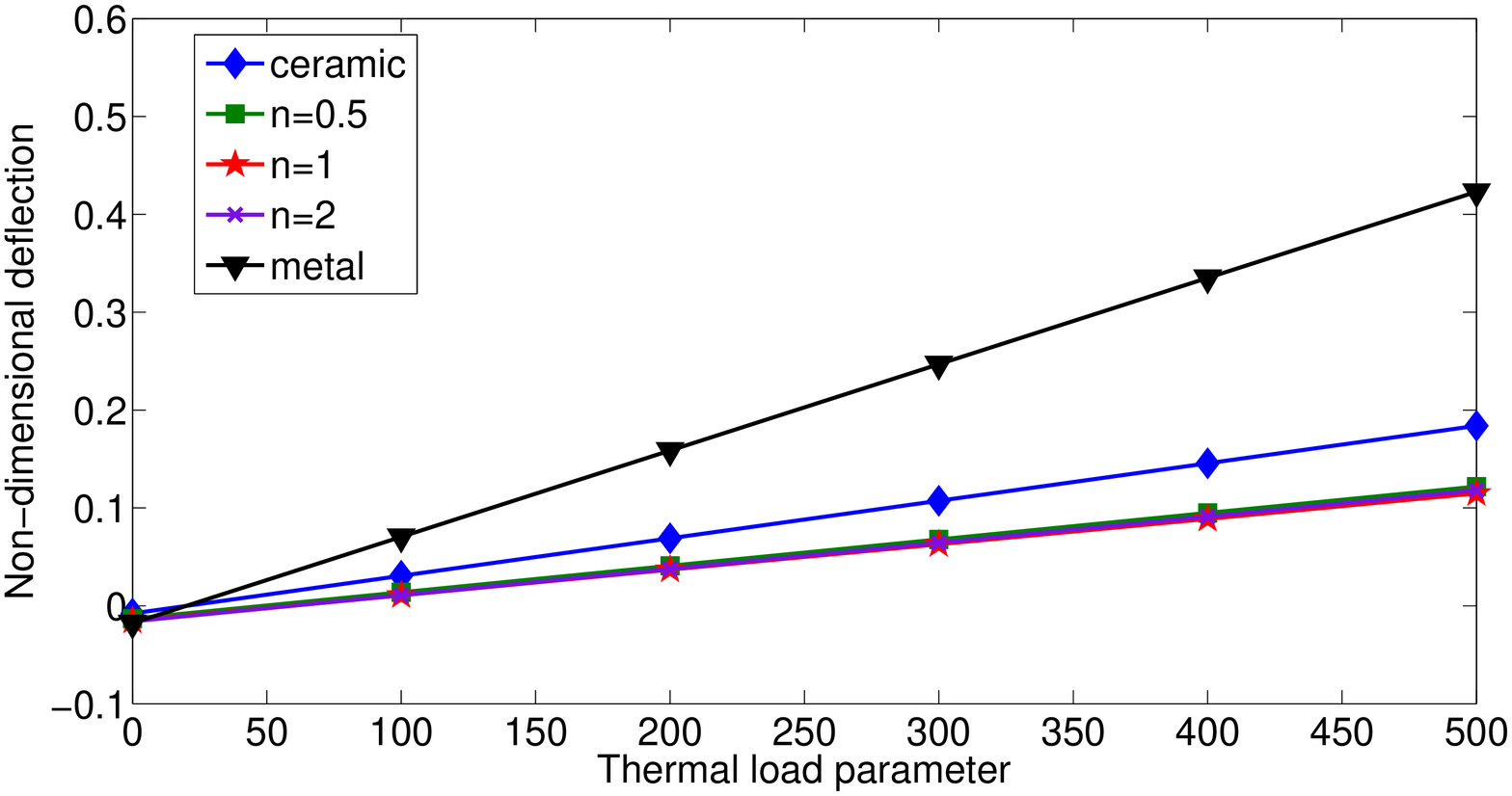}
\caption{The normalized center deflection $\overline{w}_c = w_c/h$ of a simply supported square Al/ZrO$_2$ FGM plate subjected to thermal load. }
\label{fig:SS_normDisp_thermal}
\end{figure}

Next, we study the behaviour of FGM plates under the thermal loading. A Al/ZrO$_2$-2 FGM plate with length $a=$ 0.2m and thickness $h=$ 0.01m is considered. The temperature at the top surface of the plate is varied from 0$^\circ$C to 500$^\circ$C, while the temperature at the bottom surface is maintained at 20$^\circ$C. The temperature of the stress free state is assumed to be at $T_o=$ 0$^\circ$C. \fref{fig:SS_normDisp_thermal} depicts the non-dimensional center deflection of the plate under the thermal load. This problem was solved by Zhao and Liew, using the element free $kp$-Ritz method~\cite{zhaoliew2009} and by Xuan \textit{et al.,}~\cite{nguyen-xuantran2012} using NS-DSG3 elements. The results obtained with our IGA method agree very well with those given in~\cite{zhaoliew2009,nguyen-xuantran2012}. From \fref{fig:SS_normDisp_thermal}, it is observed that the metal plate gives the maximum deflection because of its high thermal conductivity. The deflection of the FGM plate with the gradient index $n=1$ is minimum. Generally, it can be seen that the deflections of the FGM plates are much lower than those of isotropic plates $(i.e., n = 0)$, which implies high temperature resistance behaviour of the functionally graded plates.

\begin{figure}[htpb]
\centering
\includegraphics[scale=0.5]{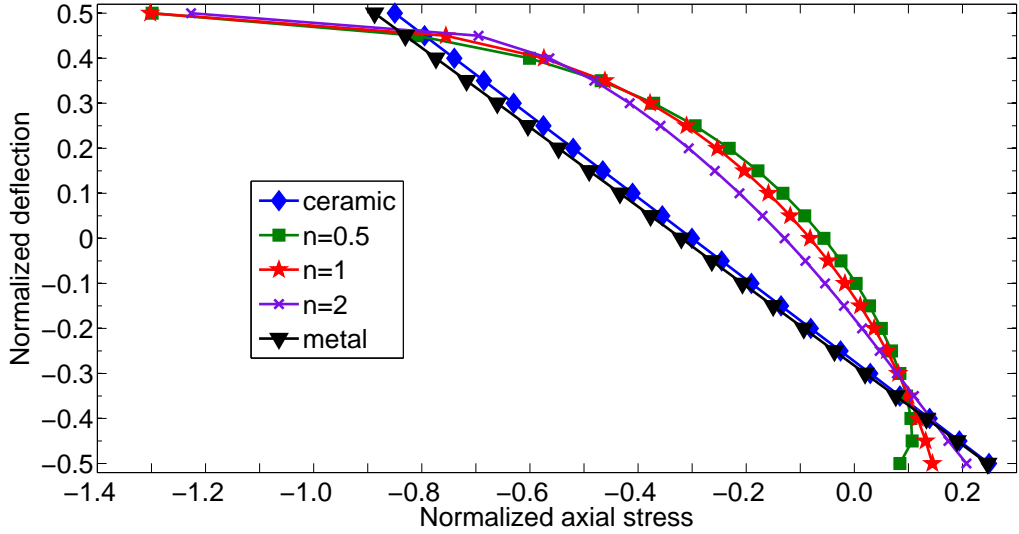}
\caption{The normalized center deflection $\overline{w}_c = w_c/h$ of a simply supported square Al/ZrO$_2$ FGM plate subjected to thermo-mechanical load. }
\label{fig:SS_normDisp_thermoMech_fgm2}
\end{figure}

Now, we investigate the FGM plate under thermo-mechanical loads. The temperature at the top surface of the plate is held at 300$^\circ$C (top surface is assumed to be rich in ceramic) and the temperature at the bottom surface (assumed to be rich in metal) is $T_m=$20$^\circ$C. \fref{fig:SS_normDisp_thermoMech_fgm2} shows the center deflection $\overline{w}_c = w_c/h$ of the plate with respect to various load parameters $\overline{p} = pa^4/(E_mh^4)$ given in the interval -14 $\le \overline{p} \le$ 0 for different values of the gradient index. It can be seen that the central deflection of the plate is completely different from the case with a purely mechanical loading  (see \fref{fig:SS_centerline_stres_fgm2}). However, similarly to the case of the FGM plate under pure mechanical loading, the center deflection of the plate linearly increases with the load. The metallic phase shows the maximum range of deflection changes and  the ceramic plate the least.

\begin{figure}[htpb]
\centering
\includegraphics[scale=0.5]{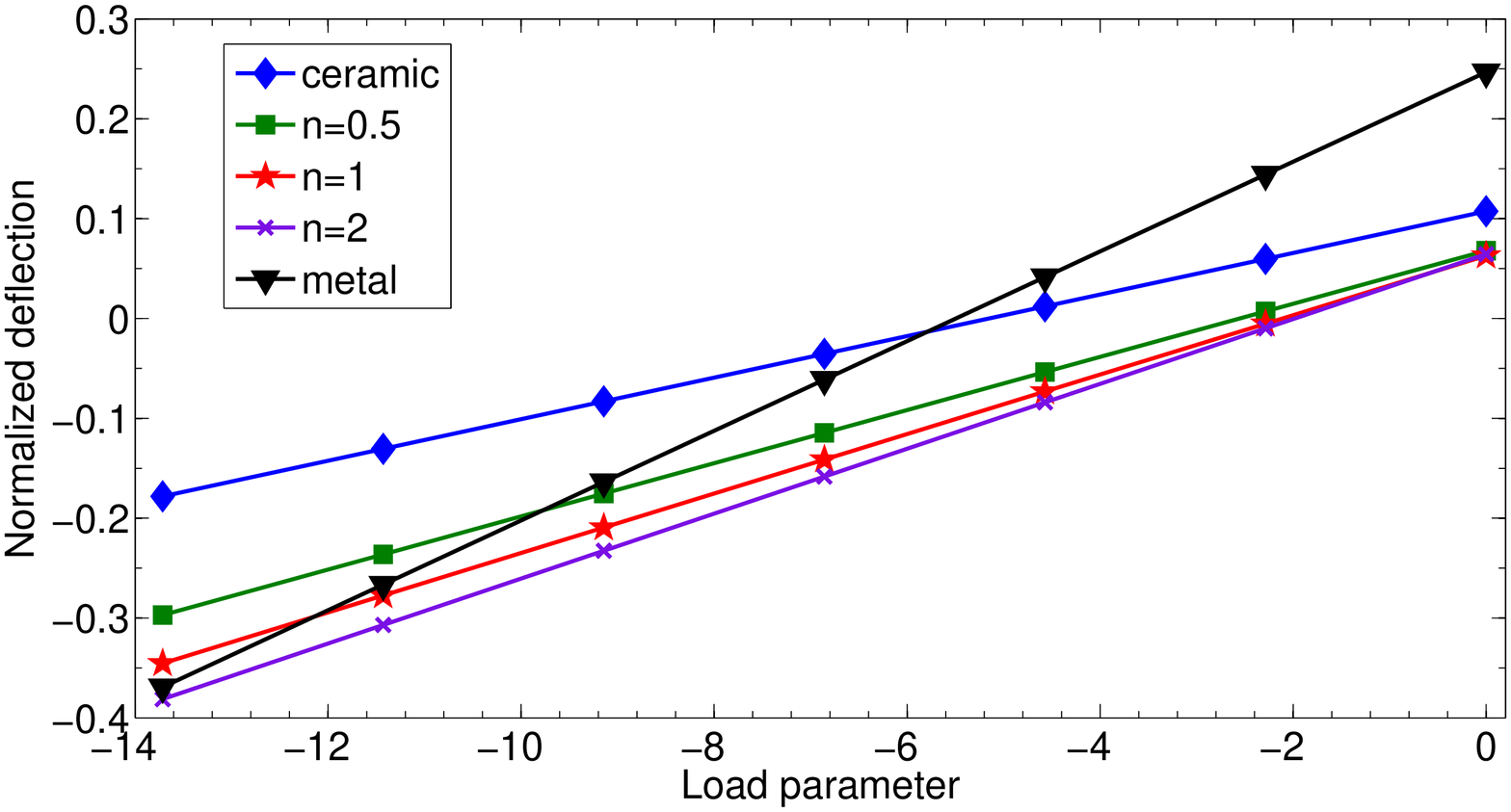}
\caption{The normalized center deflection $\overline{w}_c = w_c/h$ of a simply supported square Al/ZrO$_2$ FGM plate subjected to thermo-mechanical load. }
\label{fig:SS_normDisp_thermal_fgm2}
\end{figure}

\fref{fig:SS_normDisp_thermoMech_fgm2} plots the distribution of axial stress through the thickness of the plate under the uniform mechanical load $p=$ -10$^6$ N/m$^2$. Comparing to \fref{fig:SS_centerline_stres_fgm2} for purely mechanical load, it is seen that the maximum compressive stress at the top surface of the plate has been obtained for the FGM plate with a gradient index $n=$ 1. Again, the metallic or ceramic plate shows the minimum tensile stress at the bottom surface. Note that the results match very well with those given in~\cite{zhaoliew2009,nguyen-xuantran2012}.

\subsection*{Skew plates}
\begin{figure}[htpb]
\centering
\includegraphics[scale=0.5]{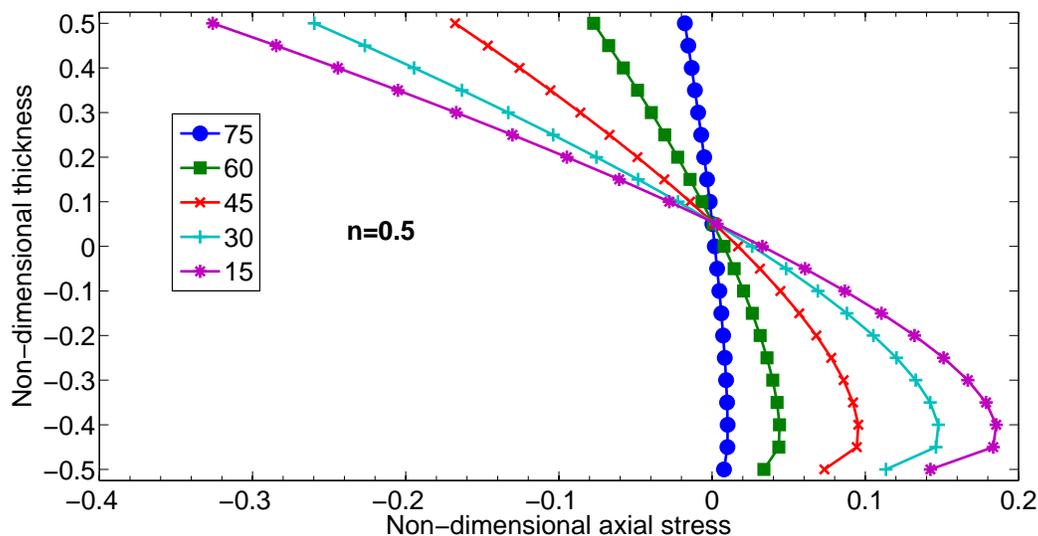}
\caption{The effect of skew angle on the axial stress $\overline{\sigma}_{xx} = \sigma_{xx} \frac{h^2}{pL^2}$ profile for simply supported FGM Al/ZrO$_2$ subjected to a mechanical load with gradient index $n=$ 0.5.}
\label{fig:skewMechLoadn05}
\end{figure}

\begin{figure}[htpb]
\centering
\includegraphics[scale=0.5]{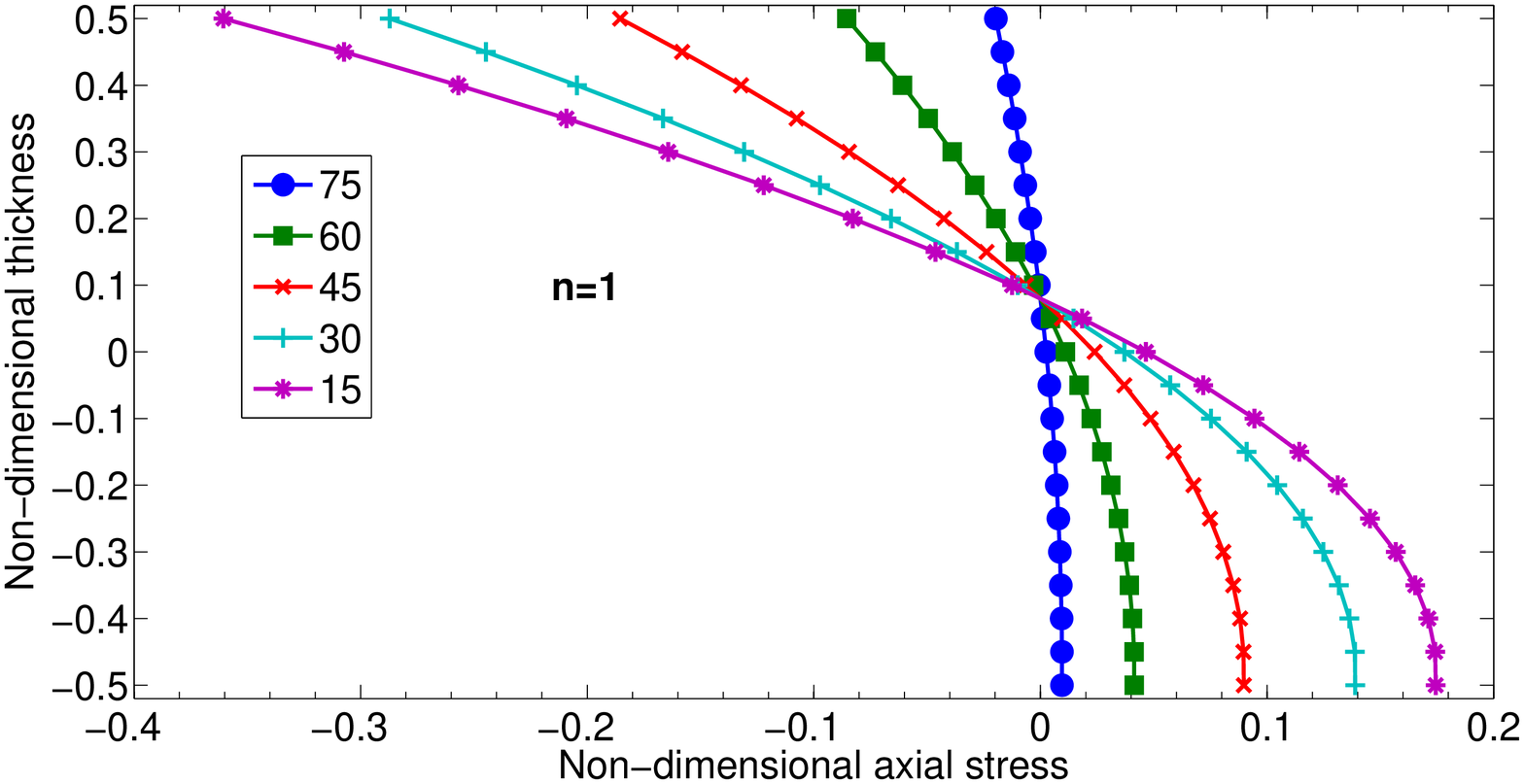}
\caption{The effect of skew angle on the axial stress $\overline{\sigma}_{xx} = \sigma_{xx} \frac{h^2}{pL^2}$ profile for simply supported FGM Al/ZrO$_2$ subjected to a mechanical load with gradient index $n=$ 1.0.}
\label{fig:skewMechLoadn1}
\end{figure}

\begin{figure}[htpb]
\centering
\includegraphics[scale=0.5]{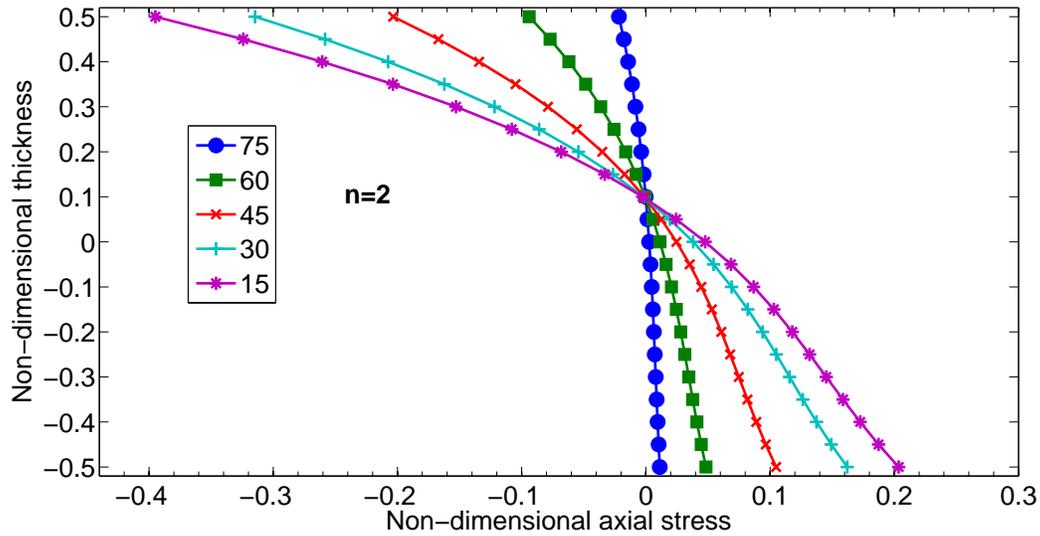}
\caption{The effect of skew angle on the axial stress $\overline{\sigma}_{xx} = \sigma_{xx} \frac{h^2}{pL^2}$ profile for simply supported FGM Al/ZrO$_2$ subjected to a mechanical load with gradient index $n=$ 2.0.}
\label{fig:skewMechLoadn2}
\end{figure}

In this example, we study the behaviour of FGM skew plates under mechanical loads. A simply supported Al/ZrO$_2$ FGM skew plate with length $a=$ 10m and thickness $h=$ 0.1m is considered. The plate is subjected to a uniform mechanical load $p=$ -10$^4$ N/m$^2$. A mesh of quadratic NURBS elements with 17 $\times$ 17 control points is used for modelling the plate. \fref{fig:skewMechLoadn05} shows the distribution of non-dimensional axial stress $\overline{\sigma}_{xx} = \sigma_{xx} h^2/(pa^2)$ through the thickness of the plate for different skew angles with gradient index $n=$ 0.5. It can be observed that the axial stresses increase as the skew angle decreases. Similar behaviour can be found for gradient indices $n=$ 1 and $n=$ 2 in \frefs{fig:skewMechLoadn1} and \ref{fig:skewMechLoadn2}, respectively. The results obtained by our isogeometric analysis are in a good agreement with those reported in \cite{leezhao2009} using the element free $kp-$ Ritz method and the results of NS-DSG3~\cite{nguyen-xuantran2012} and ES-DSG3~\cite{nguyen-xuantran2011}.

\subsection{Free flexural vibrations}
In this section, the free flexural vibration characteristics of FGM plates are studied numerically. In all cases, we present the non-dimensionalized free flexural frequency defined as, unless otherwise stated:

\begin{equation}
\overline{\omega} = \omega h \sqrt{ \frac{\rho_c}{E_c}}
\end{equation}

where $\omega$ is the natural frequency, $\rho_c, E_c$ are the mass density and Young's modulus of the ceramic phase. Before proceeding with a detailed study on the effect of gradient index on the natural frequencies, the formulation developed herein is validated against available analytical/numerical solutions in the literature.

\subsubsection*{Square plates}
\begin{figure}[htpb]
\centering
\includegraphics[scale=0.5]{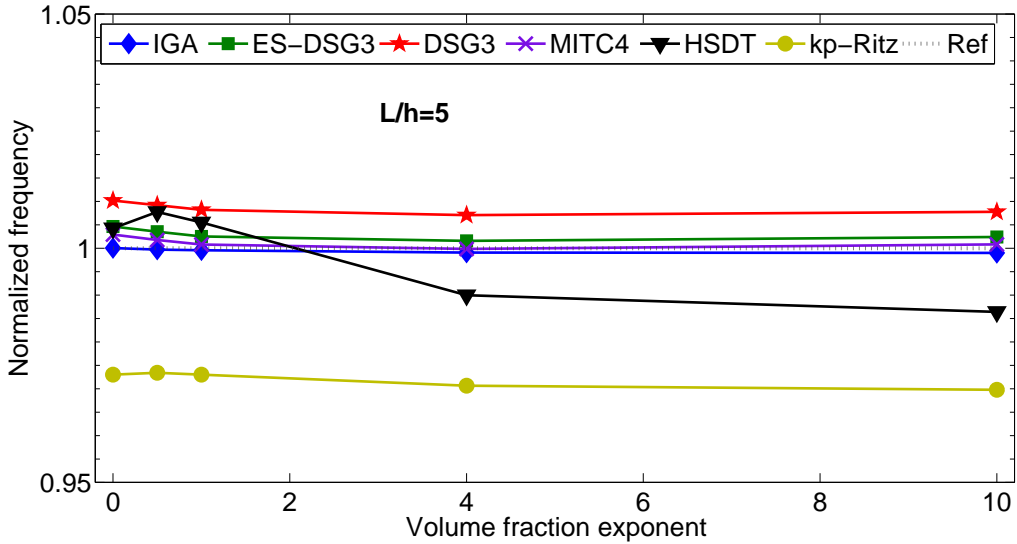}
\caption{The normalized first frequency $\omega^\ast = \omega/\omega_{\rm{ref}}$ versus volume fraction exponents for SSSS Al/Al$_2$O$_3$ FGM plate with $a/h=$ 5. The frequency is normalized with the analytical solution.}
\label{fig:normfreqah5}
\end{figure}

\begin{figure}[htpb]
\centering
\includegraphics[scale=0.5]{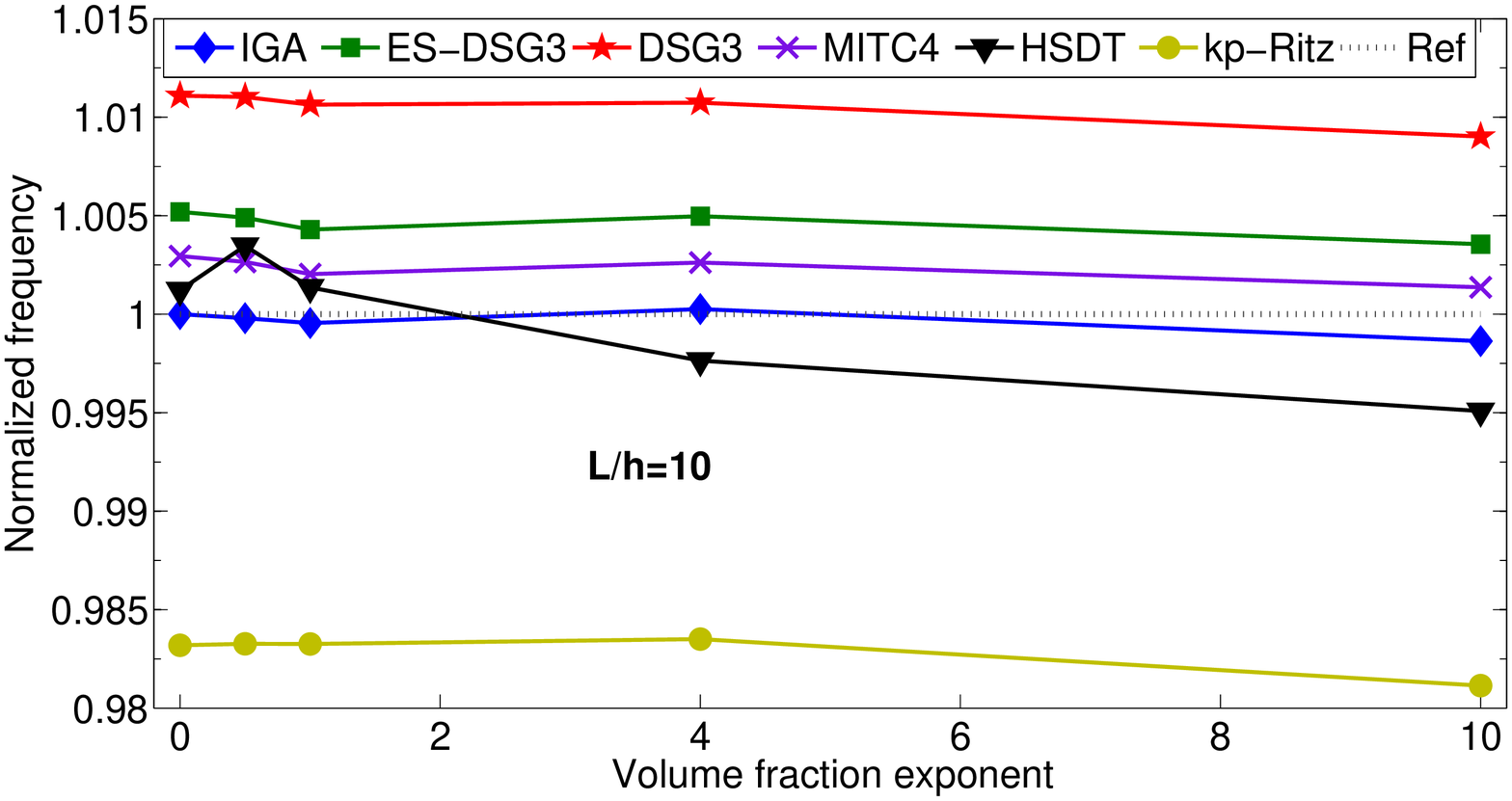}
\caption{The normalized first frequency $\omega^\ast = \omega/\omega_{\rm{ref}}$ versus volume fraction exponents for SSSS Al/Al$_2$O$_3$ FGM plate with $a/h=$ 10. The frequency is normalized with the analytical solution.}
\label{fig:normfreqah10}
\end{figure}

\begin{figure}[htpb]
\centering
\includegraphics[scale=0.5]{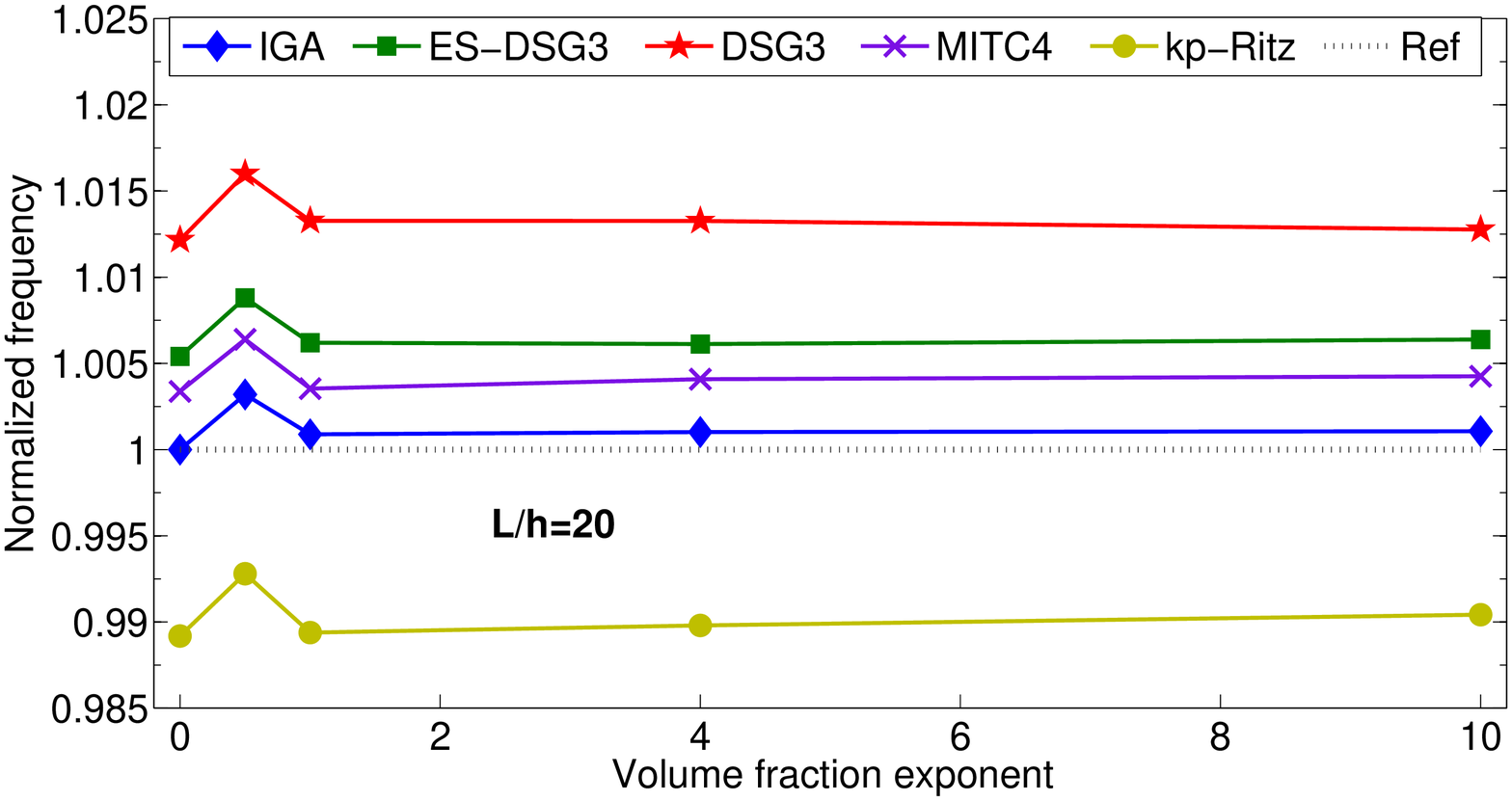}
\caption{The normalized first frequency $\omega^\ast = \omega/\omega_{\rm{ref}}$ versus volume fraction exponents for SSSS Al/Al$_2$O$_3$ FGM plate with $a/h=$ 20. The frequency is normalized with the analytical solution.}
\label{fig:normfreqah20}
\end{figure}

A  simply supported Al/Al$_2$O$_3$ FGM square plate with various length-to-thickness ratio is considered. The plate is modelled employing quadratic, cubic and quartic NURBS elements with meshes of 8, 14 and 20 control points per side. The results obtained from the isogeometric analysis for the first normalized frequency parameter $\overline{\omega}$ are presented in Tables \ref{tab:SSfgm1_freq1} and \ref{tab:SSfgm1_freq2} for different plate aspect ratios $(a/h=5,10,20)$ and compared with the results available in the literature~\cite{matsunaga2008,nguyen-xuantran2011,zhaolee2009,hashemifadaee2011}. The first normalized frequency $(= \overline{\omega}/\omega_{analytical})$ is shown in \frefs{fig:normfreqah5} - \fref{fig:normfreqah20} for various $a/h$ ratios $(a/h=5,10,20)$. The results of IGA shown in these figures are computed using quadratic NURBS elements with 14 control points per side. It is seen that the results from IGA are superior to compared to other methods, irrespective of the material gradient index.


\begin{figure}[htpb]
\centering
\subfigure[]{\includegraphics[scale=0.5]{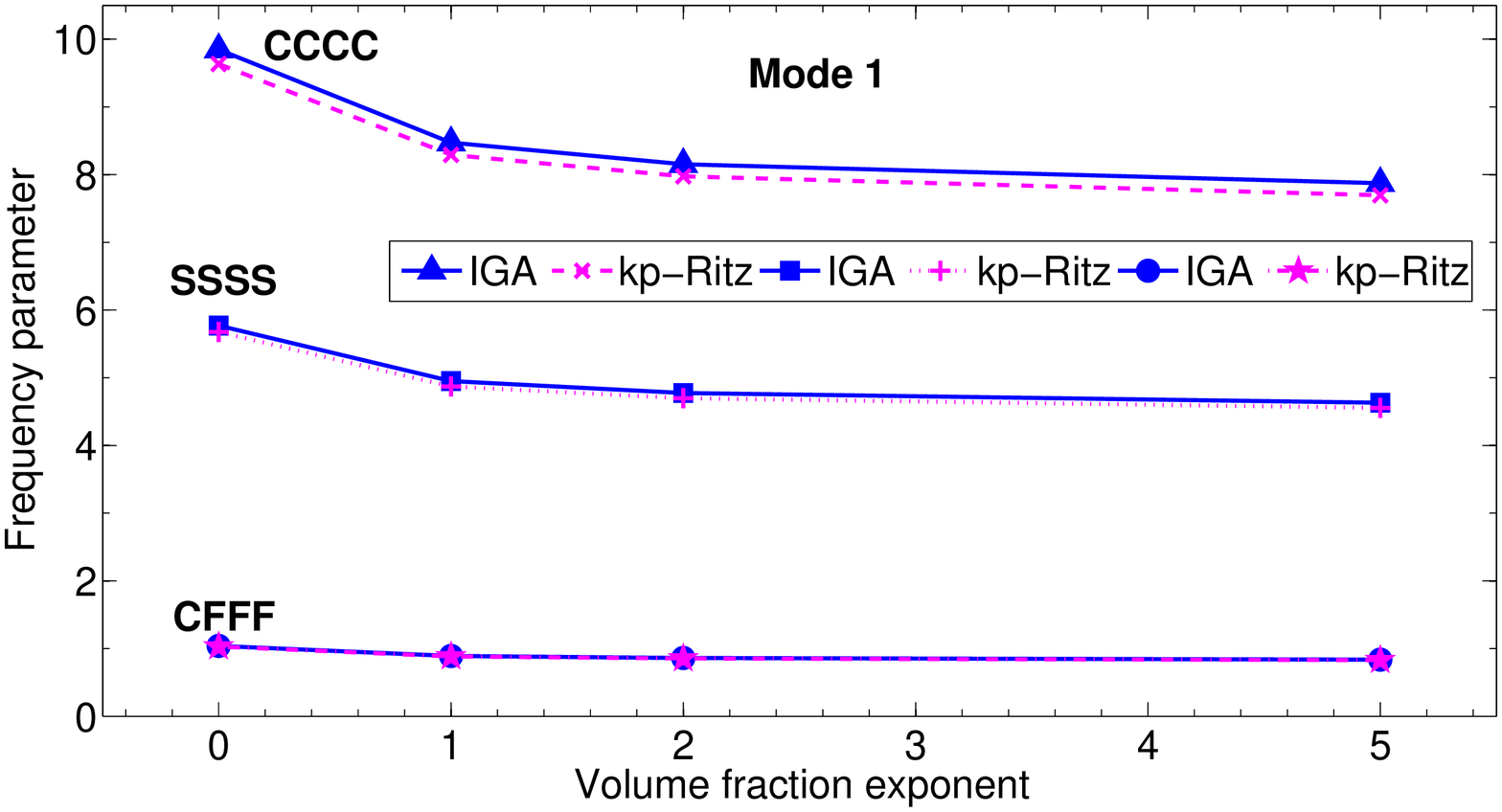}}
\subfigure[]{\includegraphics[scale=0.5]{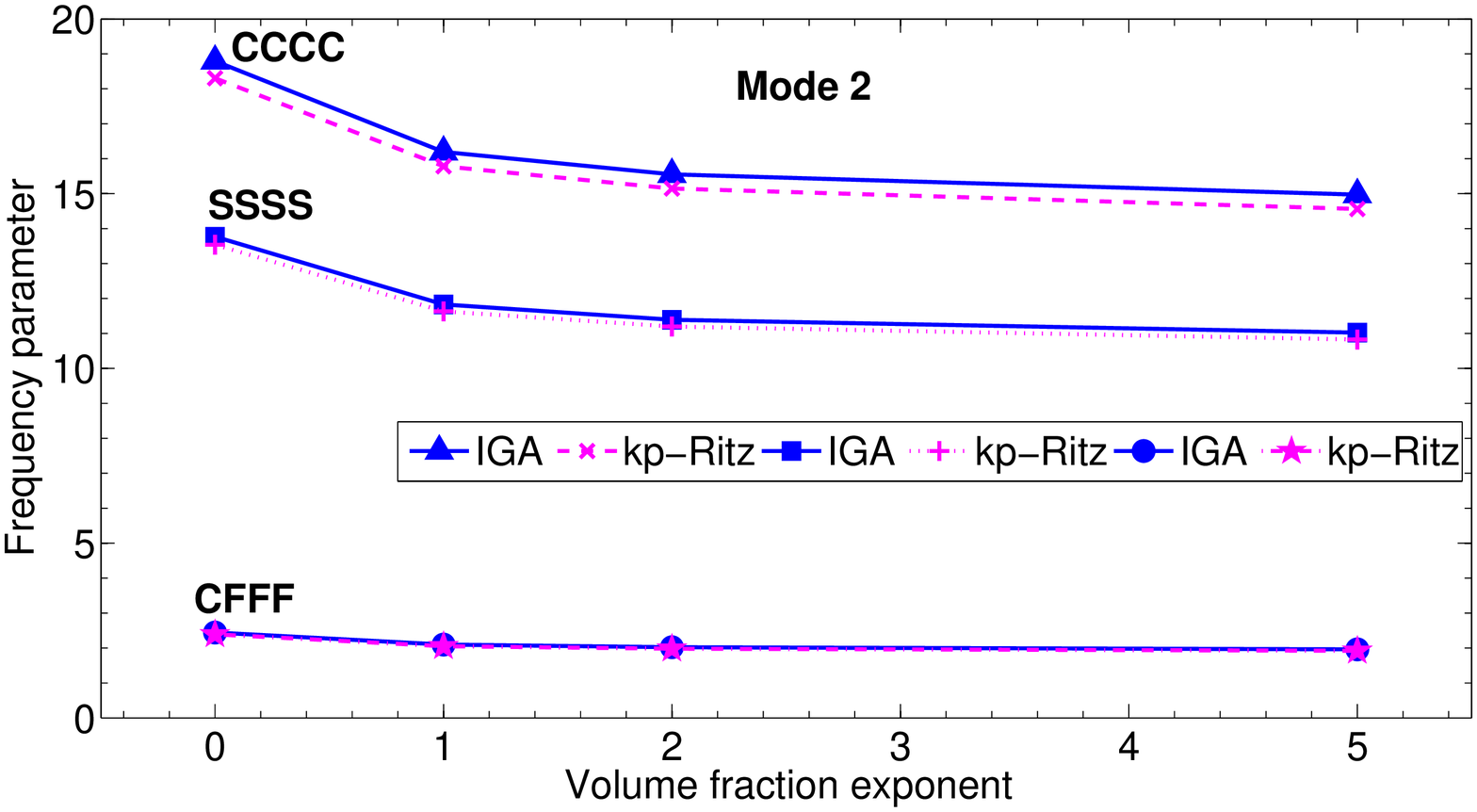}}
\caption{The non-dimensional frequency $\overline{\omega} = \frac{\omega a^2}{h} \sqrt{ \frac{\rho_c}{E_o} }$ versus various volume fraction exponent of Al/ZrO$_2$ FGM plate with different boundary conditions: (a) Mode 1 and (b) Mode 2.}
\label{fig:first2modes}
\end{figure}

\begin{figure}[htpb]
\centering
\subfigure[]{\includegraphics[scale=0.5]{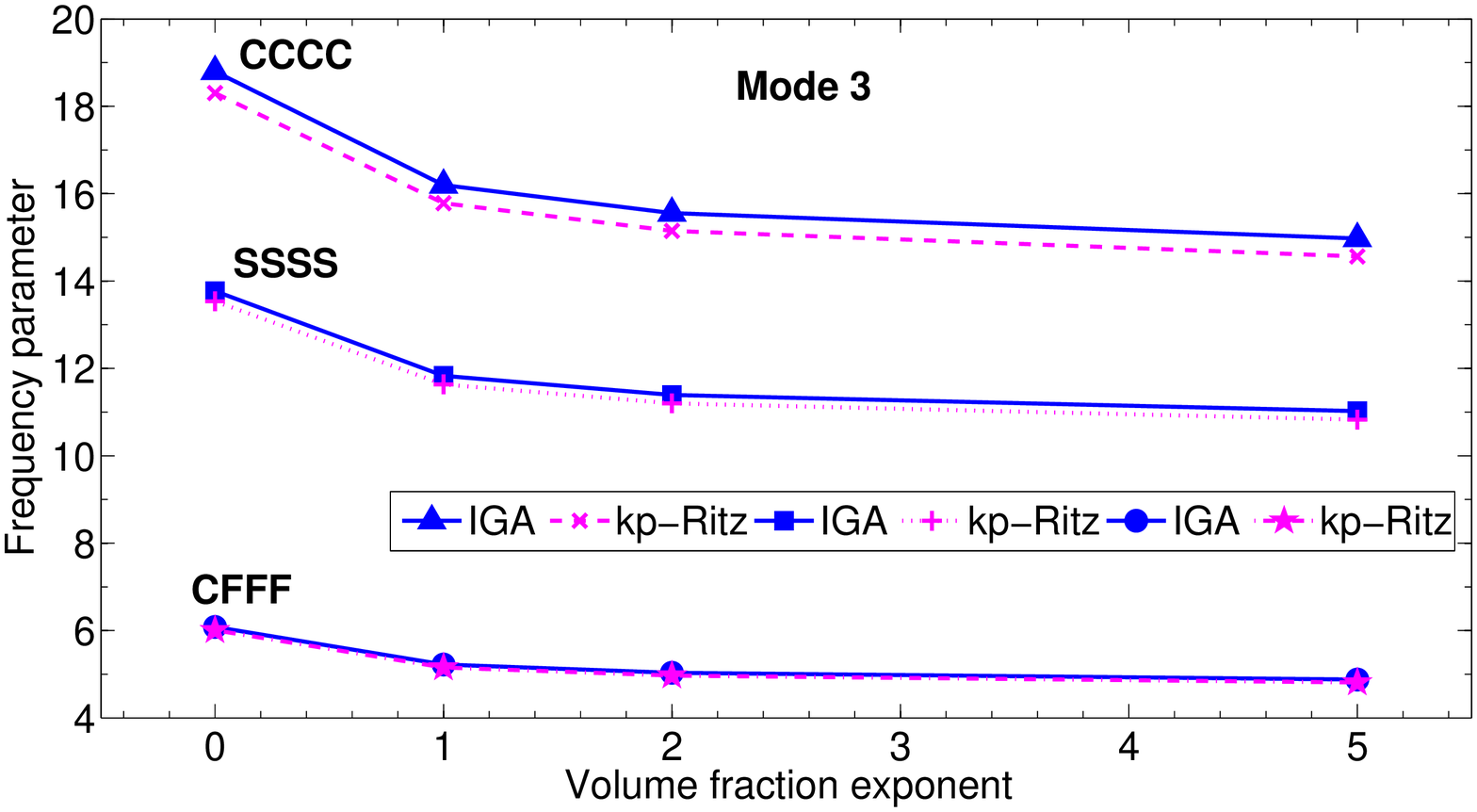}}
\subfigure[]{\includegraphics[scale=0.5]{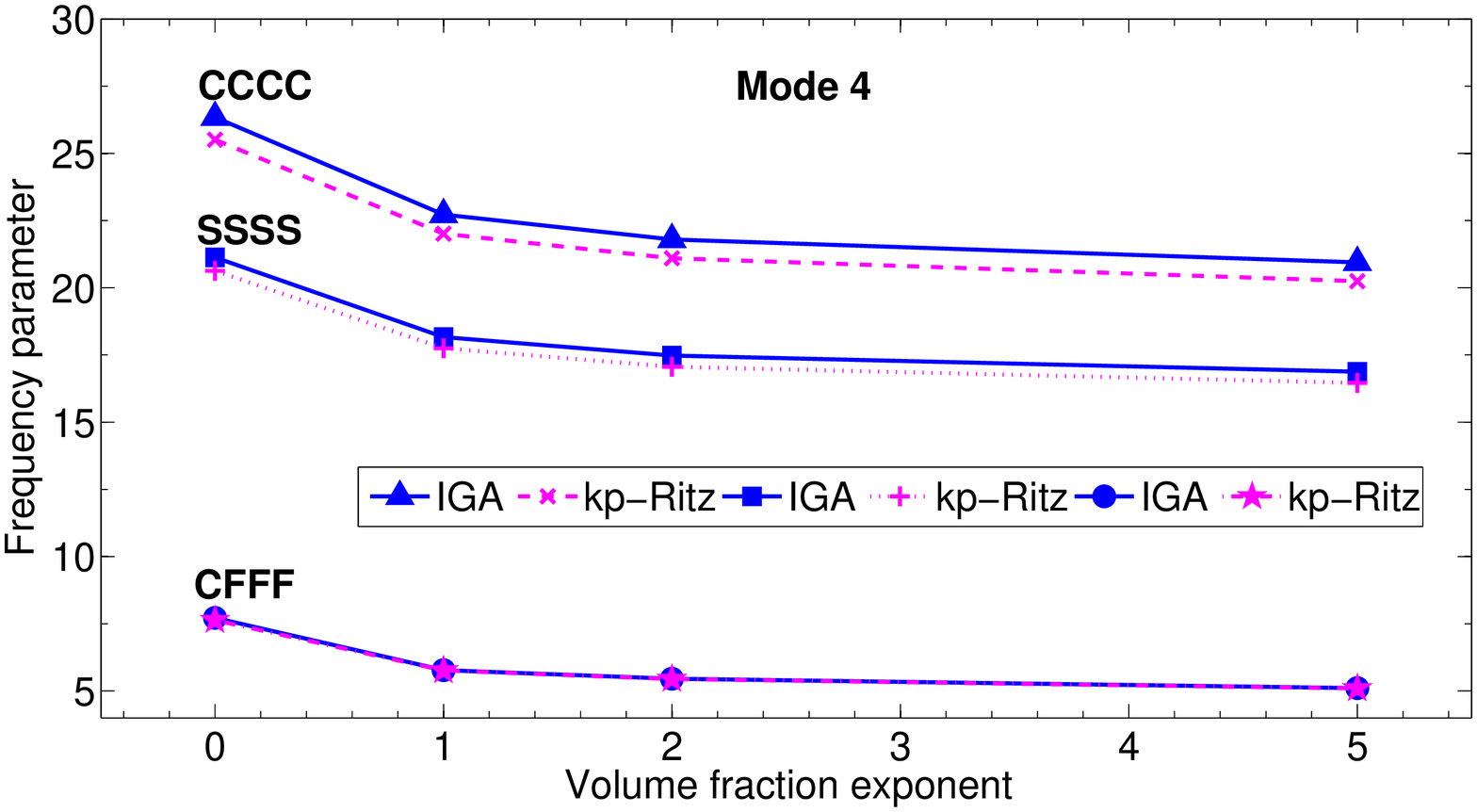}}
\caption{The non-dimensional frequency $\overline{\omega} = \frac{\omega a^2}{h} \sqrt{ \frac{\rho_c}{E_o} }$ versus various volume fraction exponent of Al/ZrO$_2$ FGM plate with different boundary conditions: (a) Mode 3 and (b) Mode 4.}
\label{fig:second2modes}
\end{figure}

\frefs{fig:first2modes} - (\ref{fig:second2modes}) shows the first four non-dimensionalized frequencies for a Al/ZrO$_2$ FGM plate with $a/h=$ 10 with various boundary conditions (CCCC, SSSS and CFFF). It can be seen that very good agreement is obtained with the results available in the literature. For CFFF boundary conditions, the difference between the IGA and the element free $kp-$Ritz is ranged from 0.3\% to 2.1\%, while it is about 2.2 - 3.3\% and 1.6 - 2.5\% for CCCC and SSSS boundary conditions, respectively. Note that all the present results agree very well with those given in~\cite{nguyen-xuantran2012,nguyen-xuantran2011}.

\begin{table}[htbp]
\centering
\renewcommand{\arraystretch}{1.5}
\caption{The first normalized frequency parameter $\overline{\omega} = \omega h \sqrt{\frac{\rho_c}{E_c}}$ for fully simply supported Al/Al$_2$O$_3$ FGM square plate for $a/h=$ 5.}
\begin{tabular}{llrrrrr}
\hline
Method & Control points & \multicolumn{5}{c}{gradient index $n$} \\
\cline{3-7}
 &  & 0 & 0.5 & 1 & 4 & 10\\
\hline
\multirow{3}{*}{IGA-Quadratic} & 8 & 0.21128 & 0.18051 & 0.16309 & 0.13962 & 0.13231 \\
& 14 & 0.21121 & 0.18045 & 0.16303 & 0.13957 & 0.13227 \\
& 20 & 0.21121 & 0.18044 & 0.16303 & 0.13957 & 0.13227 \\
\cline{3-7}
\multirow{3}{*}{IGA-Cubic} & 8 & 0.21121 & 0.18044 & 0.16303 & 0.13957 & 0.13227 \\
& 14 & 0.21121 & 0.18044 & 0.16303 & 0.13957 & 0.13227 \\
& 20 & 0.21121 & 0.18044 & 0.16303 & 0.13957 & 0.13227 \\
\cline{3-7}
\multirow{3}{*}{IGA-Quartic} & 8 & 0.21121 & 0.18044 & 0.16303 & 0.13957 & 0.13227 \\
& 14 & 0.21121 & 0.18044 & 0.16303 & 0.13957 & 0.13227 \\
& 20 & 0.21121 & 0.18044 & 0.16303 & 0.13957 & 0.13227 \\
\cline{3-7}
ES-DSG3 (20$\times$20)~\cite{nguyen-xuantran2011} & & 0.21218 & 0.18114 & 0.16351 & 0.13992 & 0.13272 \\
DSG3 (16$\times$16)~\cite{nguyen-xuantran2011} & & 0.21335 & 0.18216 & 0.16444 & 0.14069 & 0.13343 \\
MITC4 (16$\times$16)~\cite{nguyen-xuantran2011} & & 0.21182 & 0.18082 & 0.16323 & 0.13968 & 0.13251 \\
HSDT~\cite{matsunaga2008} & & 0.21210 & 0.18190 & 0.16400 & 0.13830 & 0.13060 \\
$kp-$Ritz~\cite{zhaolee2009} & & 0.20550 & 0.17570 & 0.15870 & 0.13560 & 0.12840 \\
Ref.~\cite{hashemifadaee2011} & & 0.21120 & 0.18050 & 0.16310 & 0.13970 & 0.13240 \\
\hline
\end{tabular}%
\label{tab:SSfgm1_freq1}%
\end{table}%

\begin{table}[htbp]
\centering
\renewcommand{\arraystretch}{1.5}
\caption{The first normalized frequency parameter $\overline{\omega} = \omega h \sqrt{\frac{\rho_c}{E_c}}$ for fully simply supported Al/Al$_2$O$_3$ FGM square plate for $a/h=$ 10, 20.}
\begin{tabular}{clrrrrr}
\hline
$a/h$ & Method & \multicolumn{5}{c}{gradient index $n$} \\
\cline{4-7}
 &  & 0 & 0.5 & 1 & 4 & 10\\
\hline
\multirow{9}{*}{10} & IGA-Quadratic (20 points)  & 0.05769 & 0.04899 & 0.4418 & 0.03821 & 0.03655 \\
& IGA-Cubic (20 points) & 0.05769 & 0.04898 & 0.04417 & 0.03821 & 0.03655 \\
& IGA-Quartic (20 points) & 0.05769 & 0.04898 & 0.04417 & 0.03821 & 0.03655 \\
& ES-DSG3 (20$\times$20)~\cite{nguyen-xuantran2011}  & 0.05800 & 0.04924 & 0.04439 & 0.03839 & 0.03973 \\
& DSG3 (16$\times$16)~\cite{nguyen-xuantran2011}  & 0.05834 & 0.04954 & 0.04467 & 0.03861 & 0.03693 \\
& MITC4 (16$\times$16)~\cite{nguyen-xuantran2011}  & 0.05787 & 0.049132 & 0.04429 & 0.03830 & 0.03665 \\
& HSDT~\cite{matsunaga2008}  & 0.05777 & 0.04917 & 0.04426 & 0.03811 & 0.03642 \\
& $kp-$Ritz~\cite{zhaolee2009}  & 0.05673 & 0.04818 & 0.04346 & 0.3757 & 0.03591 \\
& Ref.~\cite{hashemifadaee2011}  & 0.05770 & 0.04900 & 0.04420 & 0.03820 & 0.03660 \\
\cline{2-7}
\multirow{9}{*}{20} & IGA-Quadratic (20 points) & 0.01480 & 0.01254 & 0.01130 & 0.00981 & 0.00944 \\
& IGA-Cubic (20 points) & 0.01480 & 0.01254 & 0.01130 & 0.00981 & 0.00944 \\
& IGA-Quartic (20 points) & 0.01480 & 0.01254 & 0.01130 & 0.00981 & 0.00944 \\
& ES-DSG3 (20$\times$20)~\cite{nguyen-xuantran2011}  & 0.01488 & 0.01261 & 0.01137 & 0.00986 & 0.00946 \\
& DSG3 (16$\times$16)~\cite{nguyen-xuantran2011} &  0.01498 & 0.012704 & 0.01145 & 0.00993 & 0.00952 \\
& MITC4 (16$\times$16)~\cite{nguyen-xuantran2011}  & 0.01485 & 0.01258 & 0.01134 & 0.00984 & 0.00944 \\
& $kp-$Ritz~\cite{zhaolee2009} &  0.01464 & 0.01241 & 0.01118 & 0.00970 & 0.00931 \\
& Ref.~\cite{hashemifadaee2011} & 0.01480 & 0.01250 & 0.01130 & 0.00980 & 0.00940 \\
\hline
\end{tabular}%
\label{tab:SSfgm1_freq2}%
\end{table}%

\subsubsection*{Skew plate}
\begin{figure}[htpb]
\centering
\subfigure[]{\includegraphics[scale=0.5]{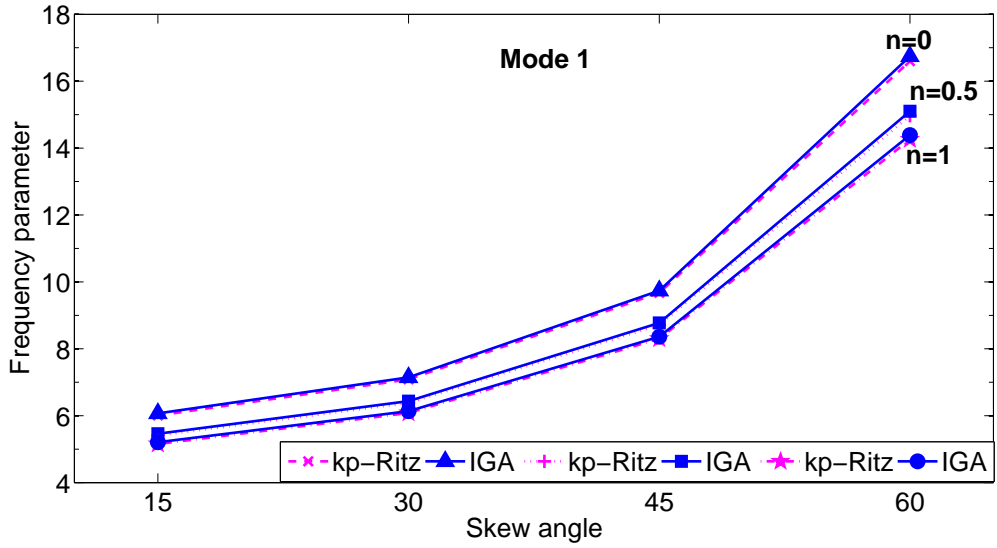}}
\subfigure[]{\includegraphics[scale=0.5]{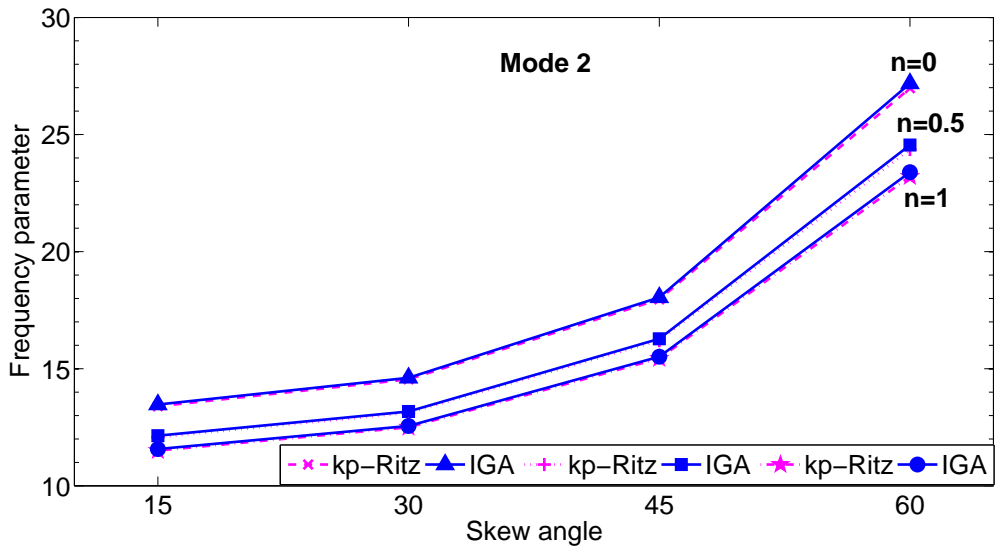}}
\caption{The non-dimensional frequency $\overline{\omega} = \omega a^2 \sqrt{\rho_c/E_c}$ for a simply supported Al/ZrO$_2$-2 skew plate versus various skew angles for different values of gradient index with $a/h=$ 10.}
\label{fig:skewFreq_ah10_SS}
\end{figure}

\begin{figure}[htpb]
\centering
\subfigure[]{\includegraphics[scale=0.5]{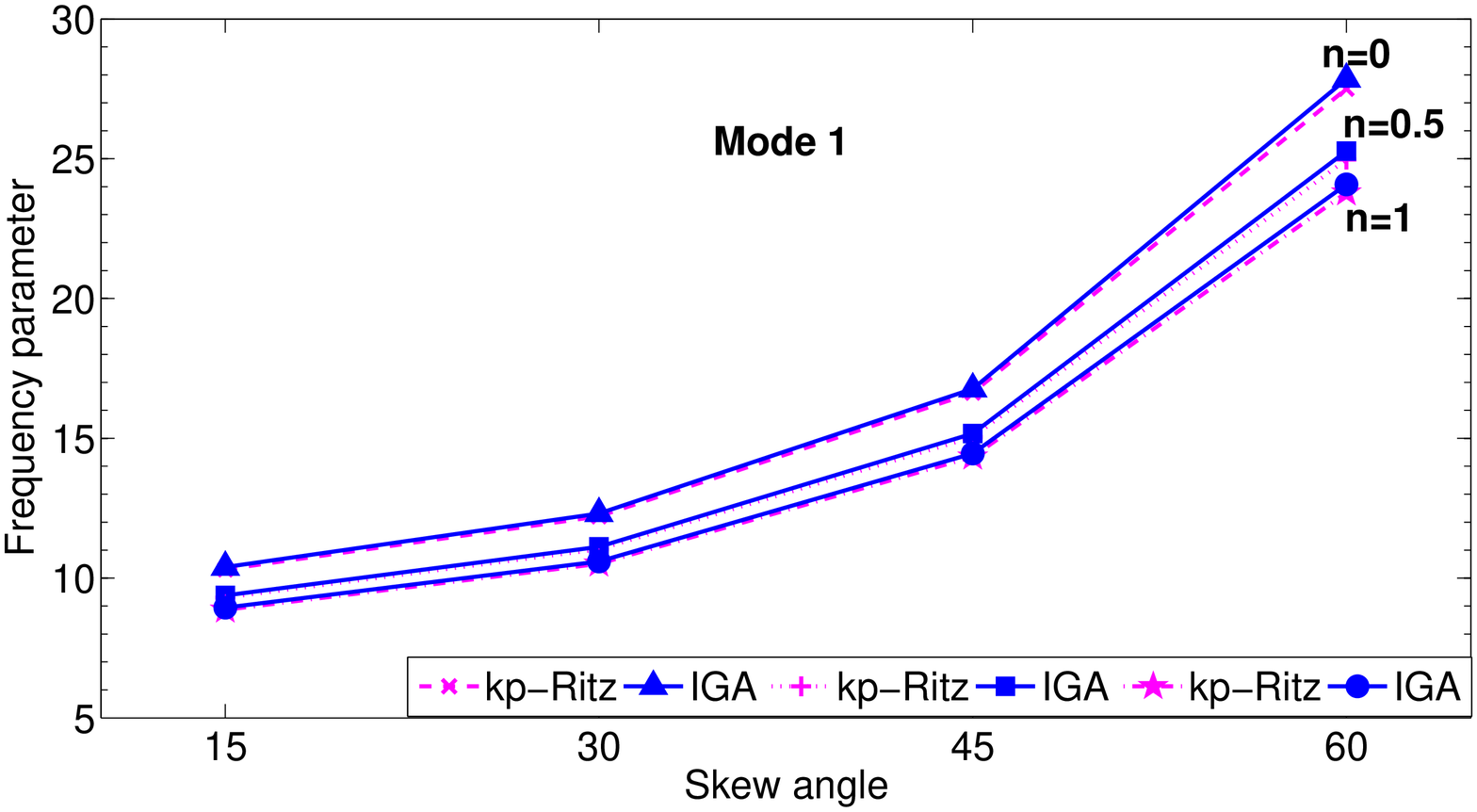}}
\subfigure[]{\includegraphics[scale=0.5]{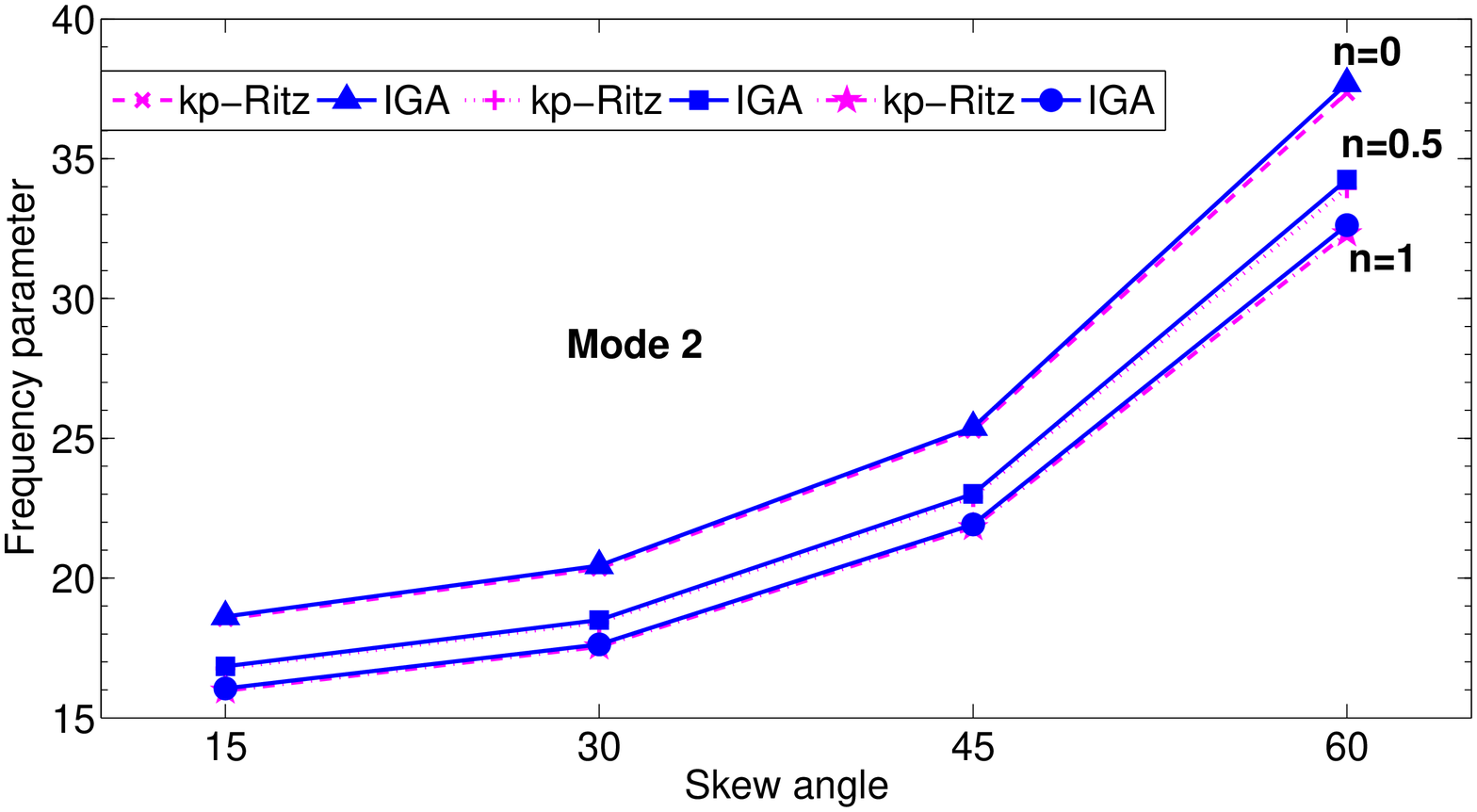}}
\caption{The non-dimensional frequency $\overline{\omega} = \omega a^2 \sqrt{\rho_c/E_c}$ for a clamped Al/ZrO$_2$-2 skew plate versus various skew angles for different values of gradient index with $a/h=$ 10.}
\label{fig:skewFreq1_ah10_CC}
\end{figure}

\begin{figure}[htpb]
\centering
\subfigure[mode 1]{\includegraphics[scale=0.25]{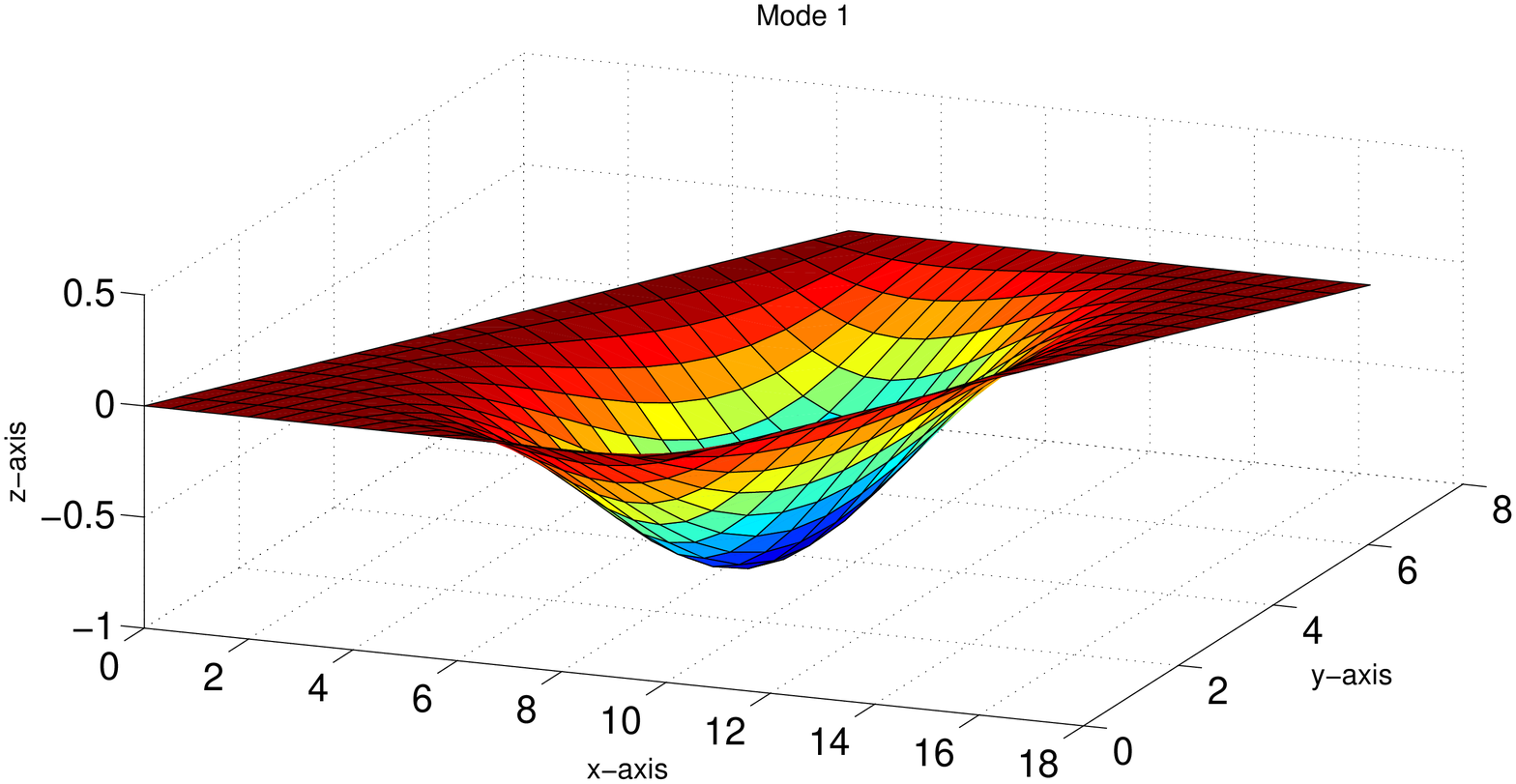}}
\subfigure[mode 2]{\includegraphics[scale=0.25]{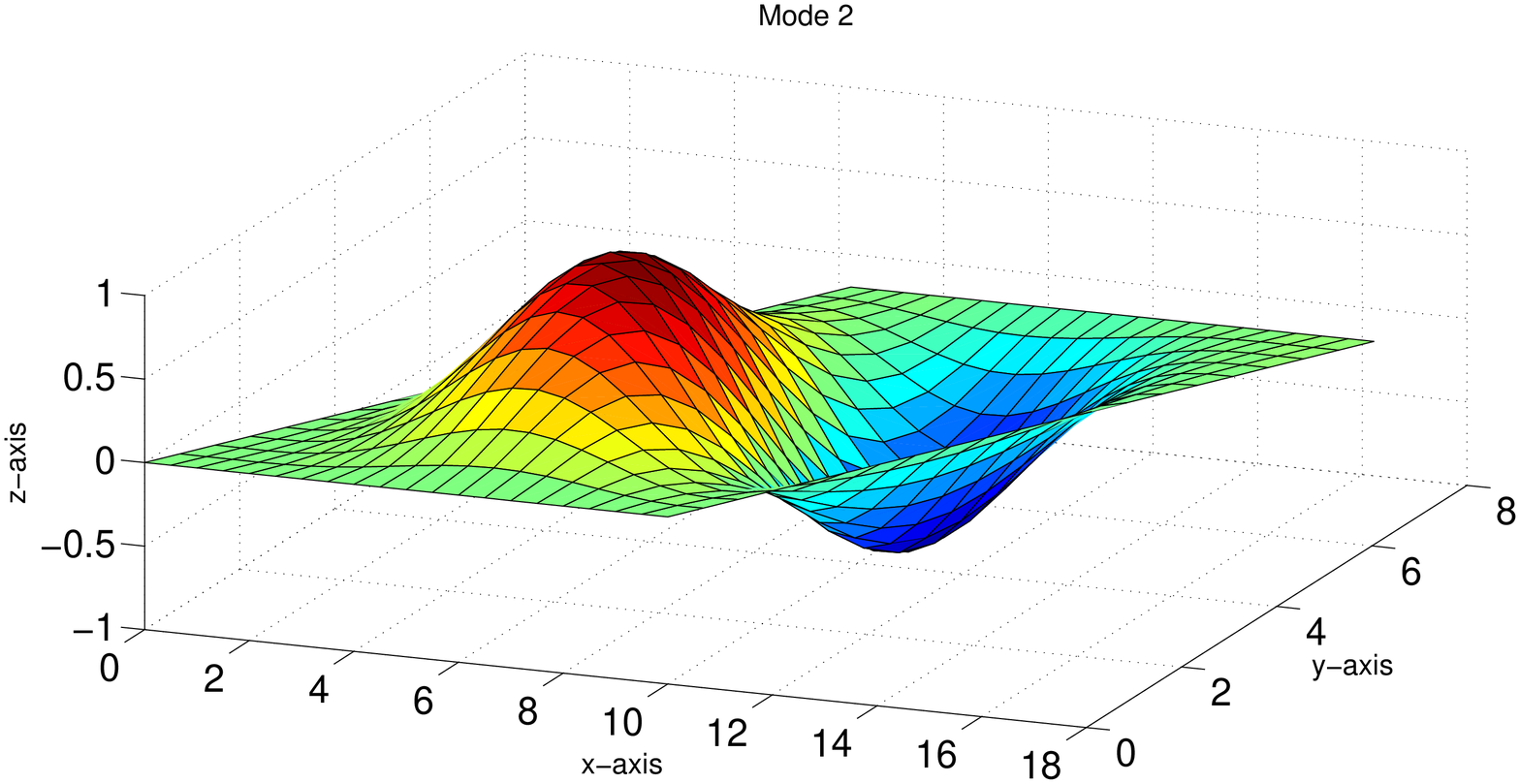}}
\subfigure[mode 3]{\includegraphics[scale=0.25]{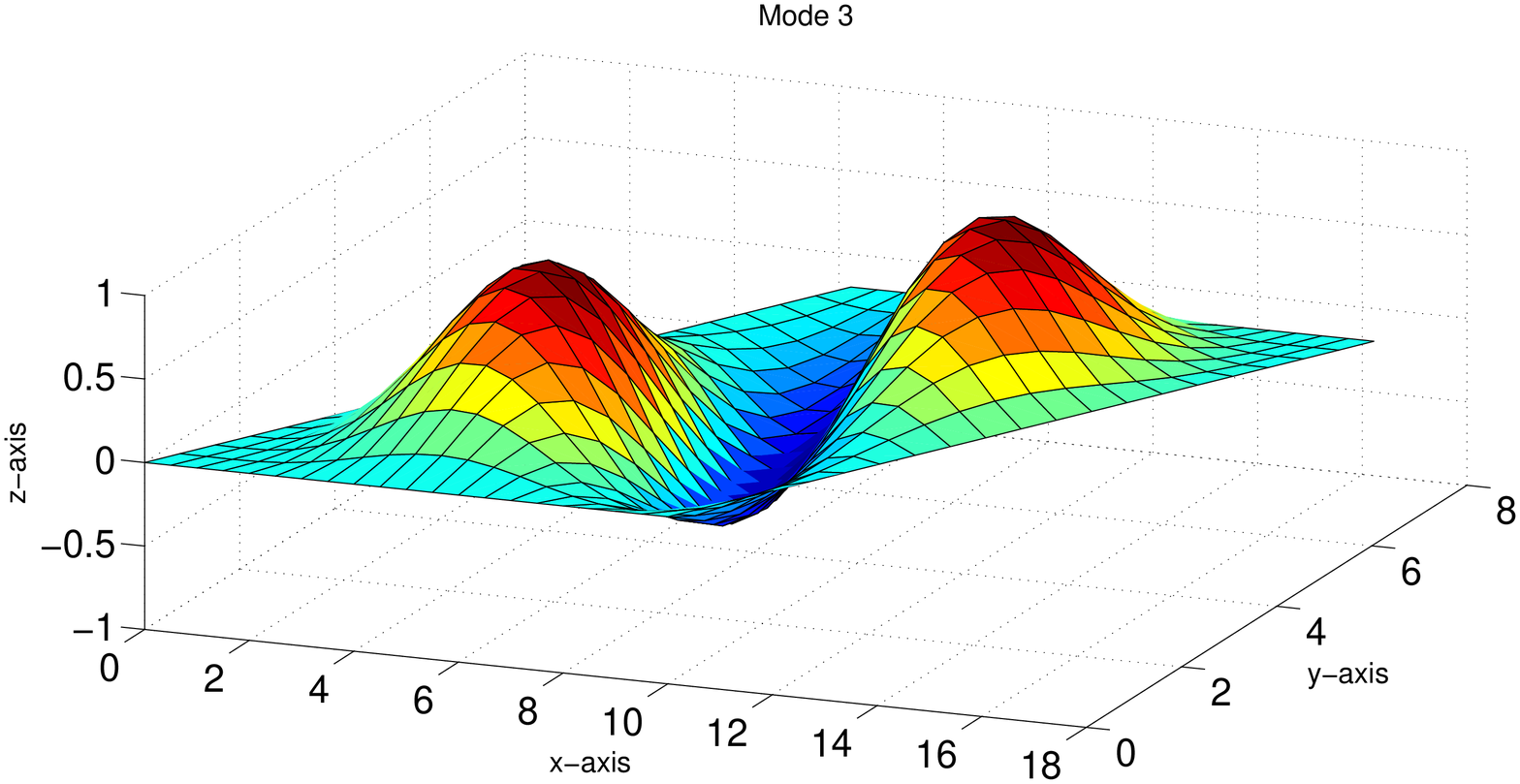}}
\subfigure[mode 4]{\includegraphics[scale=0.25]{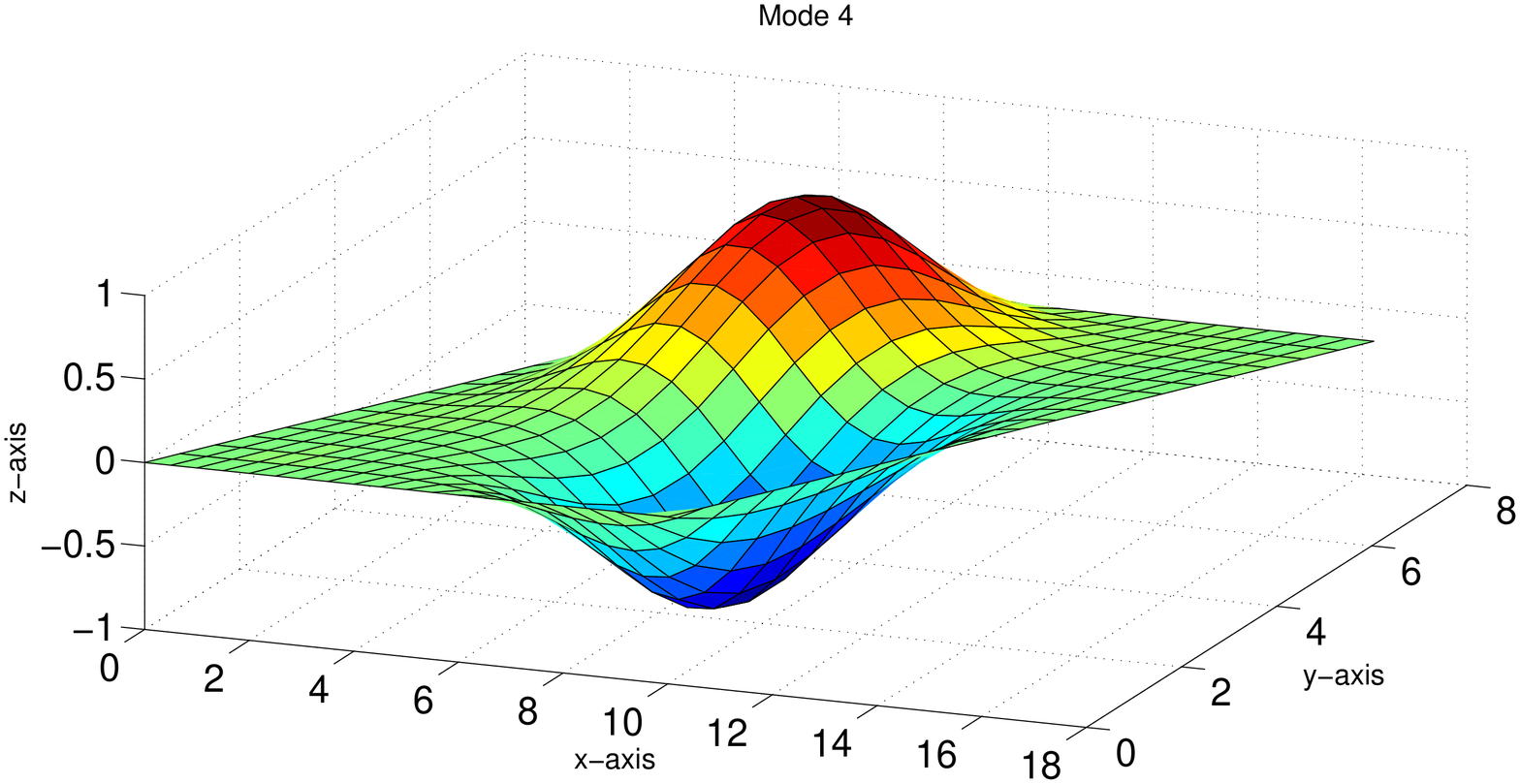}}
\subfigure[mode 5]{\includegraphics[scale=0.25]{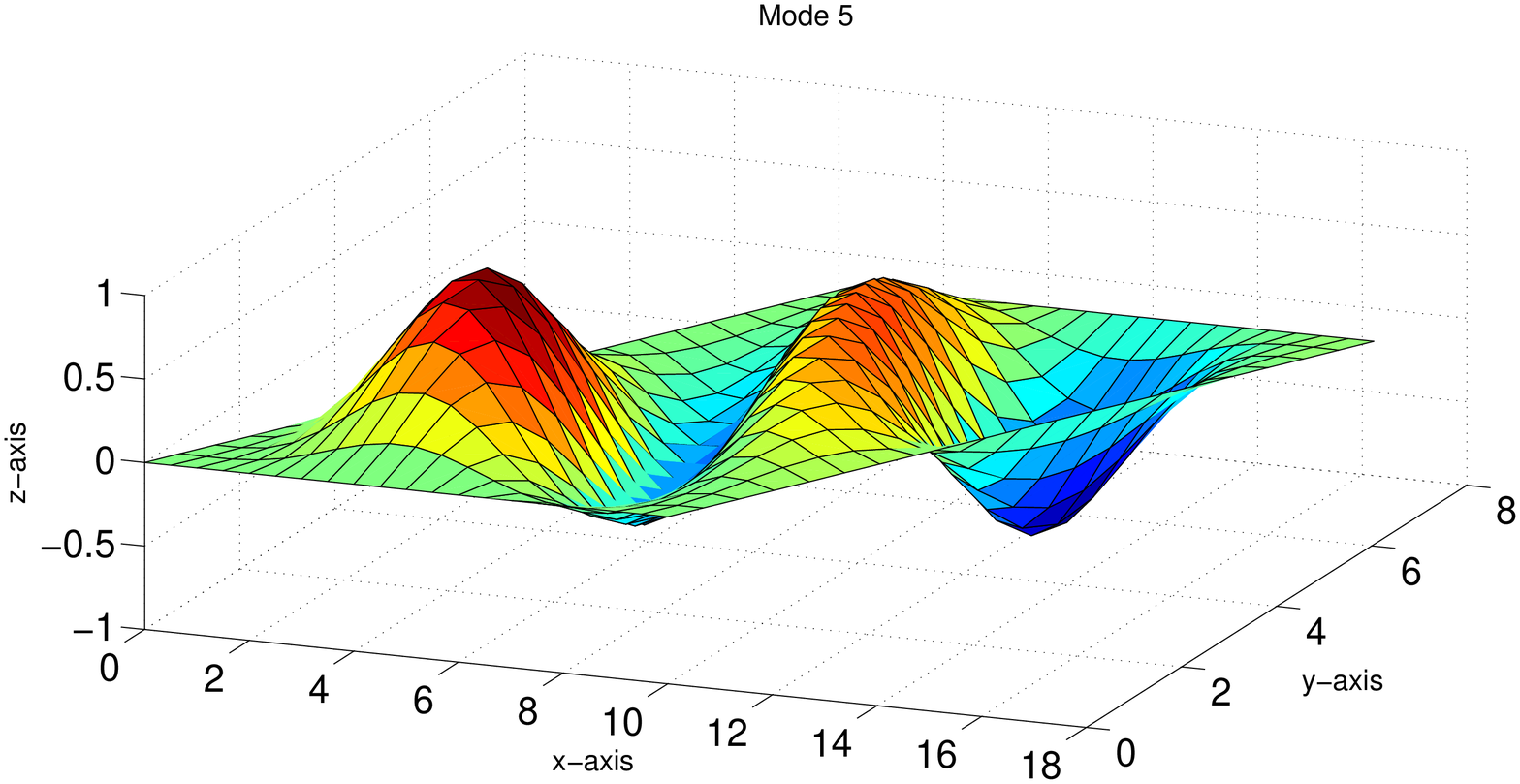}}
\subfigure[mode 6]{\includegraphics[scale=0.25]{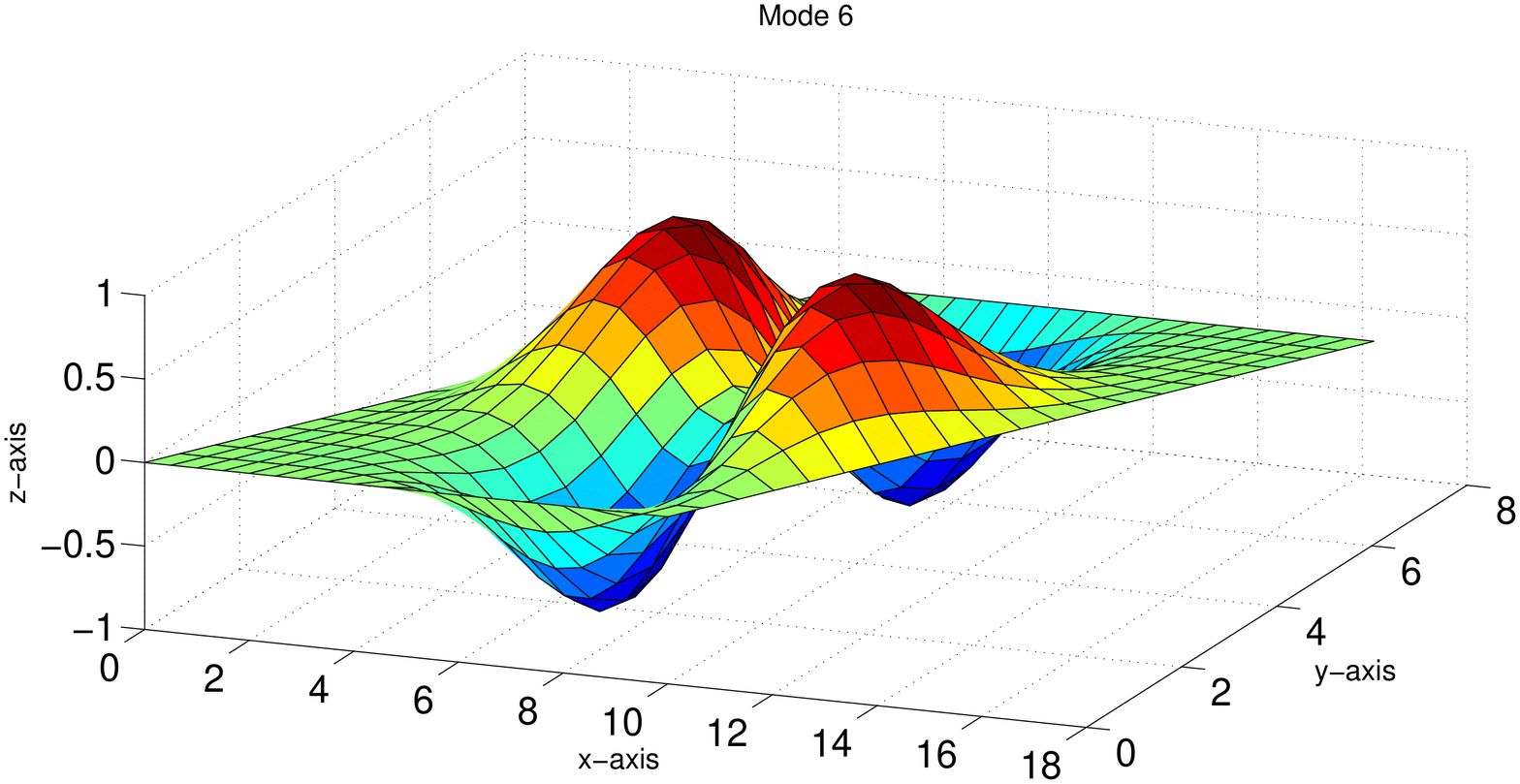}}
\subfigure[mode 7]{\includegraphics[scale=0.25]{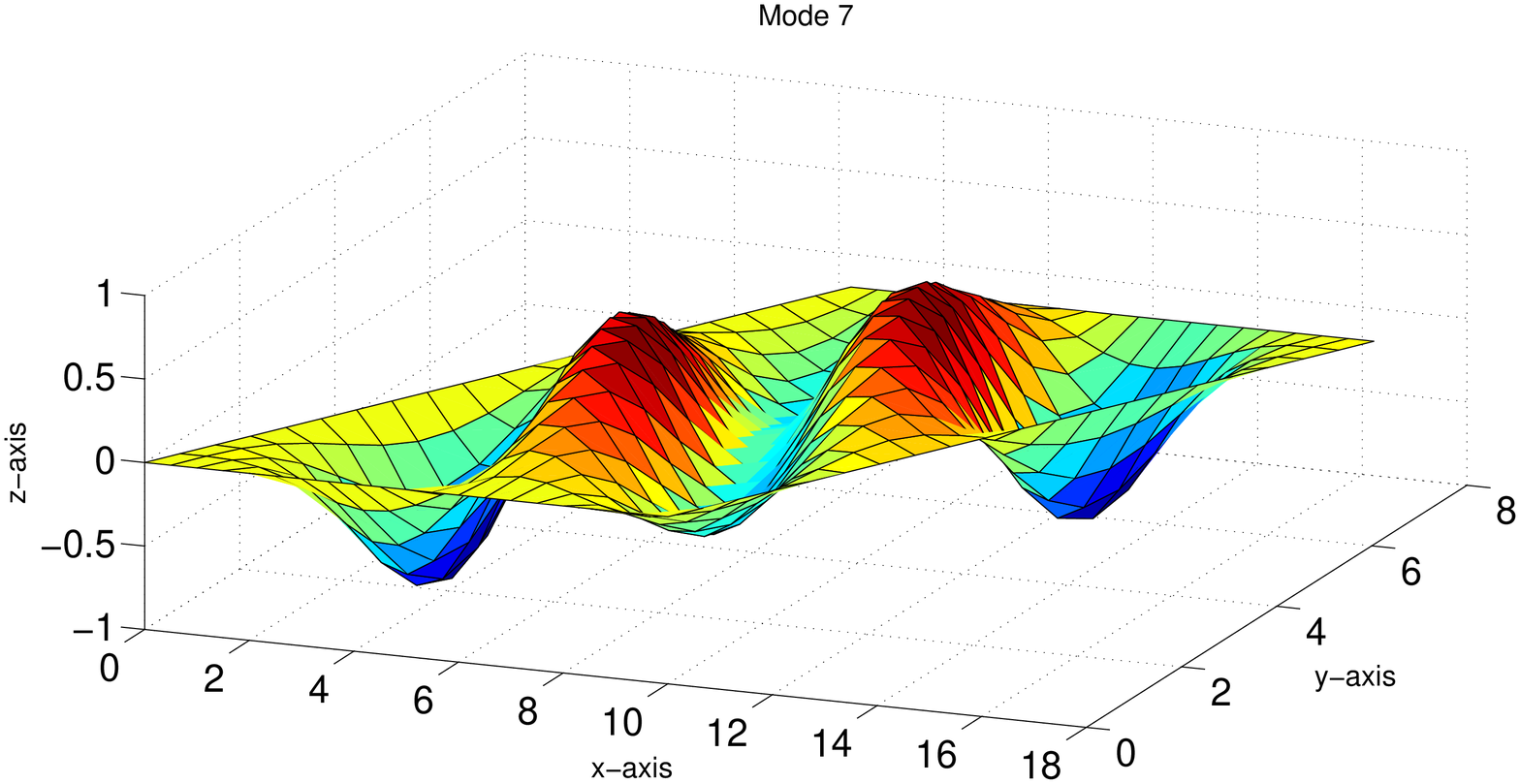}}
\subfigure[mode 8]{\includegraphics[scale=0.25]{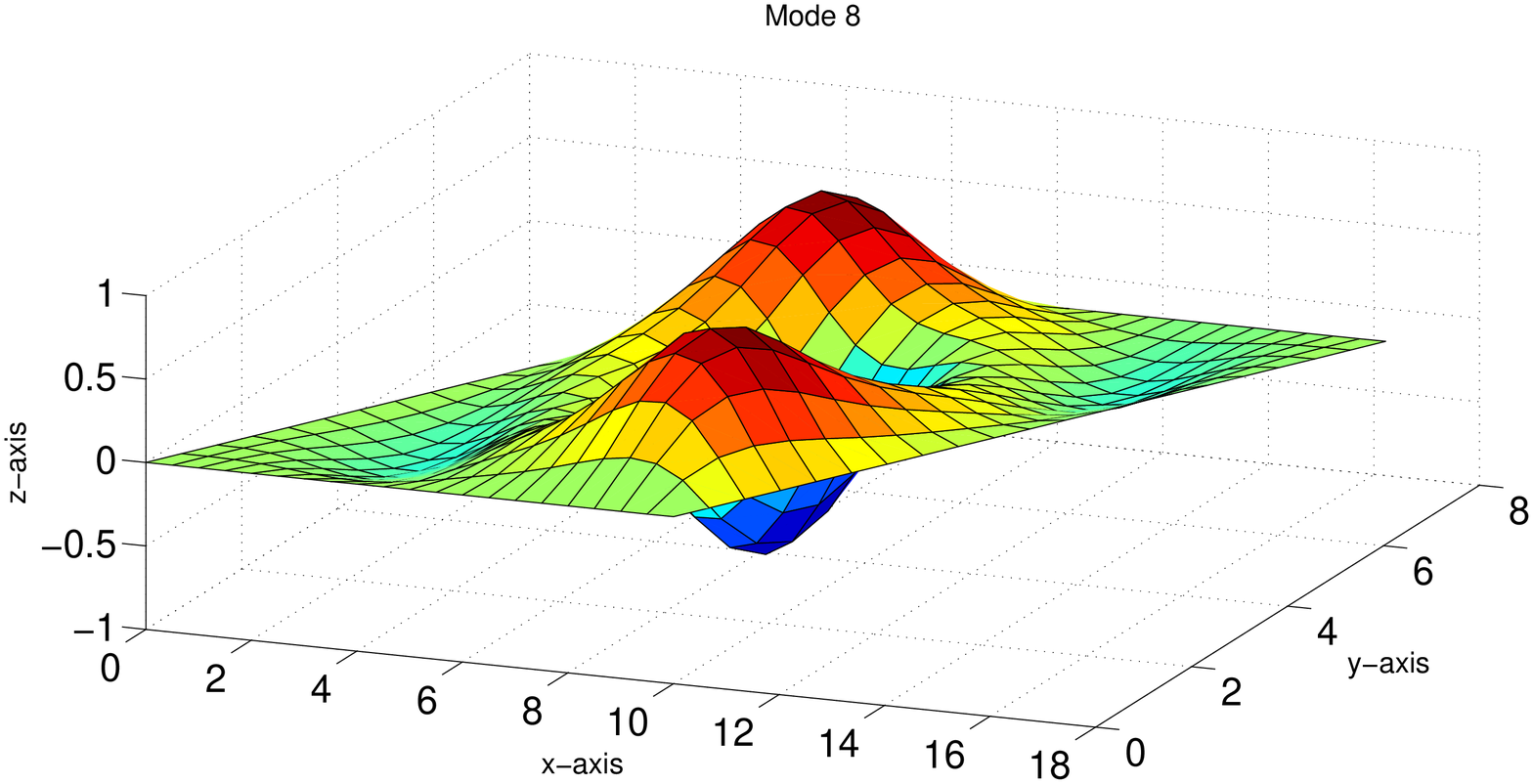}}
\caption{First eight mode shapes for a clamped Al/Zro$_2$ skew plate with $a/h=$ 10, gradient index $n=$ 0.5 and $\psi=$ 45$^\circ$.}
\label{fig:skewFreq_ah10_CC}
\end{figure}

Next, the effects of the skewness of the plate on the free flexural vibration of the FGM plate is studied. A Al/ZrO$_2$-2 skew FGM plate with length-to-thickness ratio $a/h=$ 10 and various skew angles are considered in this example. The skew plate is modelled using quadratic NURBS elements with 17 $\times$ 17 control points. The first two non-dimensionalized frequencies are shown in \frefs{fig:skewFreq_ah10_SS} - (\ref{fig:skewFreq1_ah10_CC}) for SSSS and CCCC boundary conditions, respectively. For comparison, the results from the element free $kp-$Ritz~\cite{zhaolee2009} are also plotted. From the figures, it is seen that the IGA gives higher frequencies than the element free $kp-$Ritz method. It is also observed that by increasing the gradient index or decreasing the skew angle, the frequency decreases. In both cases, the decrease in the natural frequency can be attributed to the stiffness degradation. In the case of gradient index, the stiffness degradation is due to increased metallic volume fraction, while the geometry of the plate is a contributing factor in decreasing the frequency when the skew angle decreases. The first eight mode shapes of fully clamped Al/ZrO$_2$-2 skew plate with skew angle $\psi=$ 45$^\circ$ and gradient index $n=$ 0.5 is plotted in \fref{fig:skewFreq_ah10_CC}. 

\begin{rmk}
In the next two sections, the present formulation is extended to study the buckling and flutter characteristics of the the FGM plates. For the whole, quadratic NUBRS functions with 17 $\times$ 17 control points are used, unless otherwise specified.
\end{rmk}

\subsection{Buckling analysis}
In this section, we present the mechanical and thermal buckling behaviour of functionally graded skew plates. The FGM plate considered here consists of aluminum and alumina (see Table \ref{tab:matprop} for material properties).


\subsubsection*{Mechanical Buckling}


The critical buckling parameters are defined for uni- and bi- axial compressive loads as:

\begin{eqnarray}
\lambda_{\rm{cru}} = \frac{N_{\rm{xxcr}}^0 b^2}{\pi^2 D_c} \nonumber \\
\lambda_{\rm{crb}} = \frac{N_{\rm{yycr}}^0 b^2}{\pi^2 D_c} 
\label{eqn:mbuckparm}
\end{eqnarray}

where, $D_c=E_ch^3/(12(1-\nu^2))$. The critical buckling loads evaluated by varying the skew angle of the plate, volume fraction index and considering mechanical loads such as uni- and biaxial compressive loads are shown in Tables \ref{table:mbuckpah100}  - \ref{table:mbuckpah10} for two different thickness ratios. The efficacy of the present formulation is demonstrated by comparing our results with those in~\cite{ganapathiprakash2006}. It can be seen that increasing the gradient index decreases the critical buckling load. It is also observed that the decrease in the critical value is significant for the material gradient index $n \le 2$ and that further increase in $n$ yields less reduction in the critical value, irrespective of the skew angle. 

\begin{table}[htpb]
\renewcommand\arraystretch{1.5}
\caption{Critical buckling parameters for a thin simply supported FGM skew plate with $a/h=$ 100 and $a/b=$ 1.}
\centering
\begin{tabular}{ccccccccc}
\hline 
Skew angle & $\lambda_{cr}$ & \multicolumn{7}{c}{Gradient index, $n$} \\
\cline{3-9}
 &  & \multicolumn{2}{c}{0} & \multicolumn{2}{c}{1} & 2 & 5 & 10\\
\cline{3-6}
 &  & Ref.~\cite{ganapathiprakash2006} & Present & Ref.~\cite{ganapathiprakash2006} & Present \\
\hline
\multirow{2}{*}{0$^\circ$} & $\lambda_{\rm{cru}}$ & 4.0010 & 3.9998 & 1.7956 & 1.8034 & 1.5320 & 1.2606 & 1.0830 \\
& $\lambda_{\rm{crb}}$ & 2.0002 & 1.9999 & 0.8980 & 0.9017 & 0.7660 & 0.6303 & 0.5415 \\
\multirow{2}{*}{15$^\circ$} & $\lambda_{\rm{cru}}$ & 4.3946 & 4.3946 & 1.9716 & 1.9716 & 1.6752 & 1.3800 & 1.1868 \\
& $\lambda_{\rm{crb}}$ & 2.1154 & 2.1154 & 0.9517 & 0.9517 & 0.8086 & 0.6652 & 0.5716 \\
\multirow{2}{*}{30$^\circ$} & $\lambda_{\rm{cru}}$ & 5.8966 & 5.8966 & 2.6496 & 2.6496 & 2.2515 & 1.8607 & 1.6032 \\
& $\lambda_{\rm{crb}}$ & 2.5365 & 2.5365 & 1.1519 & 1.1519 & 0.9788 & 0.8044 & 0.6905 \\
\multirow{2}{*}{45$^\circ$} & $\lambda_{\rm{cru}}$ & 10.1031 & 10.1031 & 4.5445 & 4.5445 & 3.8625 & 3.2234 & 2.7964 \\
& $\lambda_{\rm{crb}}$ & 3.6399 & 3.6399 & 1.6863 & 1.6863 & 1.4330 & 1.1774 & 1.0103 \\
\hline
\end{tabular}
\label{table:mbuckpah100}
\end{table}

\begin{table}[htpb]
\renewcommand\arraystretch{1.5}
\caption{Critical buckling parameters for a thick simply supported FGM skew plate with $a/h=$ 10 and $a/b=$ 1.}
\centering
\begin{tabular}{ccccccccc}
\hline 
Skew angle & $\lambda_{cr}$ & \multicolumn{7}{c}{Gradient index, $n$} \\
\cline{3-9}
 &  & \multicolumn{2}{c}{0} & \multicolumn{2}{c}{1} & 2 & 5 & 10\\
\cline{3-6}
 &  & Ref.~\cite{ganapathiprakash2006} & Present & Ref.~\cite{ganapathiprakash2006} & Present \\
\hline
\multirow{2}{*}{0$^\circ$} & $\lambda_{\rm{cru}}$ & 3.7374 & 3.7307 & 1.6892 & 1.6793 & 1.4198 & 1.1632 & 0.9999 \\
& $\lambda_{\rm{crb}}$ & 1.8686 & 1.8654 & 0.8449 & 0.8397 & 0.7099 & 0.5816 & 0.4999 \\
\multirow{2}{*}{15$^\circ$} & $\lambda_{\rm{cru}}$ & 4.0791 & 4.0791 & 1.8458 & 1.8458 & 1.5616 & 1.2810 & 1.1021 \\
& $\lambda_{\rm{crb}}$ & 1.9660 & 1.9660 & 0.8923 & 0.8923 & 0.7550 & 0.6184 & 0.5315 \\
\multirow{2}{*}{30$^\circ$} & $\lambda_{\rm{cru}}$ & 5.3571 & 5.3571 & 2.4298 & 2.4298 & 2.0533 & 1.6886 & 1.4565 \\
& $\lambda_{\rm{crb}}$ & 2.3226 & 2.3226 & 1.0659 & 1.0659 & 0.9011 & 0.7367 & 0.6326 \\
\multirow{2}{*}{45$^\circ$} & $\lambda_{\rm{cru}}$ & 8.5261 & 8.5261 & 3.8835 & 3.8835 & 3.2679 & 2.7046 & 2.3521 \\
& $\lambda_{\rm{crb}}$ & 3.1962 & 3.1962 & 1.5030 & 1.5030 & 1.2680 & 1.0335 & 0.8871 \\
\hline
\end{tabular}
\label{table:mbuckpah10}
\end{table}

\subsubsection*{Thermal Buckling}
The temperature rise of $T_m=$ 5$^\circ$C in the metal-rich surface of the plate is assumed in the present study. In addition to nonlinear temperature distribution across the plate thickness, the linear case is also considered in the present analysis by truncating the higher order terms in \Eref{eqn:heatconducres}. The plate is of uniform thickness and simply supported on all four edges. The critical buckling temperature difference $\Delta T_{cr}$ using two values of the aspect ratio $a/b=$ 1 and 2 with $a/h=10$ and for various skew angles is given in Table \ref{table:tbuckpah10}. It can been seen that the results from the present formulation are in good agreement with the results available in the literature. The decrease in the critical buckling load with the material gradient index $n$ is attributed to the stiffness degradation due to the increase in the metallic volume fraction. The thermal stability of the plate increases with the skew angle of the plate and the same behavior is observed for other values of gradient index $n$. It can also be seen that the nonlinear temperature variation through the thickness yields higher critical values compared to the linear distribution case.

\begin{table}[htpb]
\renewcommand\arraystretch{1.5}
\caption{Critical buckling temperature $\Delta T_{\rm{cr}}$ for a thin simply supported FGM skew plate with $a/h=$ 100 and $a/b=$ 1 under linear and nonlinear temperature rise through the thickness of the plate.}
\centering
\begin{tabular}{cclccccc}
\hline
$a/b$ & Skew angle & Temperature rise & \multicolumn{5}{c}{Gradient index, $n$} \\
\cline{4-8}
 &  &  & \multicolumn{2}{c}{0} & & &\\
\cline{4-5}
 &  &  & Ref.~\cite{ganapathiprakash2006a} & Present & 0.5 & 1 & 5 \\
\hline
\multirow{6}{*}{1} & \multirow{2}{*}{0$^\circ$} & Linear & 24.1951 & 24.1912 & 9.3787 & 5.5207 & 3.8987 \\
& & Nonlinear & 24.1951 & 24.1912 & 12.3629 & 7.6615 & 4.8740 \\
& \multirow{2}{*}{30$^\circ$} & Linear & 33.9558 & 33.9503 & 14.9115 & 9.7737 & 7.4681 \\
& & Nonlinear & 33.9558 & 33.9503 & 19.6600 & 13.558 & 9.3399 \\
& \multirow{2}{*}{60$^\circ$} & Linear & 123.0974 & 123.1172 & 65.4519 & 48.6271 & 40.0647 \\
& & Nonlinear & 123.0974 & 123.1172 & 86.2949 & 67.4989 & 50.1064 \\
\cline{3-8}
\multirow{3}{*}{2} & 0$^\circ$ & Nonlinear & 75.4278 & 75.4475 & 50.6564 & 38.6525 & 28.3067 \\
& 30$^\circ$ & Nonlinear & 100.9349 & 100.9512 & 69.7225 & 54.0838 & 39.9718 \\
& 60$^\circ$ & Nonlinear & 304.6912 & 304.6421 & 222.0921 & 177.4280 & 133.1456 \\
\hline
\end{tabular}
\label{table:tbuckpah10}
\end{table}

\subsection{Flutter analyses}
In this section, the present formulation is extended to analyse the flutter characteristics of functionally graded material plates. Both simply supported and clamped boundary conditions are considered in this study and the flow direction is assumed to be at right angles to the plate. Only square plate is considered and the results are presented only for $a/h=100$. It should be noted that the present formulation is not limited to this alone. In all cases, we present the non dimensionalized critical aerodynamic pressure, $\lambda_{cr}$ and critical frequency $\omega_{cr}$ as, unless specified otherwise:

\begin{eqnarray}
\Omega_{cr} = \omega_{cr} a^2 \sqrt{ \frac{\rho_c h}{D_c}} \nonumber \\
\lambda_{cr} = \lambda_{cr} \frac{a^3}{D_c}
\label{eqn:nondimfreq}
\end{eqnarray}

where $D_c = {E_c h^3 \over 12(1-\nu_c^2)}$ is the bending rigidity of the plate, $E_c, \nu_c$ are the Young's modulus and Poisson's ratio of the ceramic material and $\rho_c$ is the mass density. In order to be consistent with the existing literature, properties of the ceramic are used for normalization.


\begin{table}[htpb]
\renewcommand\arraystretch{1.5}
\caption{Comparison of critical aerodynamic pressure and coalescence frequency for an isotropic plate with various boundary conditions $(a/b=1, a/h=100, \nu=0.3, \theta^\prime = 0)$.}
\centering
\begin{tabular}{ccrr}
\hline 
Reference & Flutter bounds & \multicolumn{2}{c}{Boundary condition} \\
\cline{3-4}
1& 1 & Simply supported & Clamped \\
\hline
\multirow{2}{*}{Ref.~\cite{Prakash2006}} & $\lambda_{cr}$ & 511.11 & 852.34 \\
& $\omega_{cr}$ & 1840.29 & 4274.32 \\
\multirow{2}{*}{Present} & $\lambda_{cr}$ & 511.92 & 854.88 \\
& $\omega_{cr}$ & 1844.80 & 4305.30\\
\hline
\end{tabular}
\label{table:isoFlutComp}
\end{table}

\begin{figure}[htpb]
\centering
\includegraphics[scale=0.5]{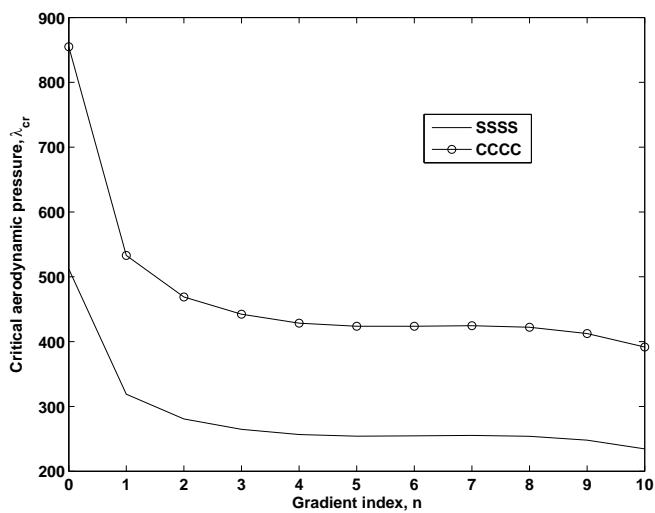}
\caption{Effect of boundary conditions on the critical aerodynamic pressure $\lambda_{cr}$}
\label{fig:flutFGMRes}
\end{figure}

Before proceeding with the detailed study, the formulation developed herein is validated against available results pertaining to the critical aerodynamic pressure and critical frequency for an isotropic plate with and without a crack. The computed critical aerodynamic pressure and the critical frequency for an isotropic square plate with various boundary conditions is given in  Table~\ref{table:isoFlutComp}. Next, the influence of boundary conditions on the flutter characteristics is studied. For this study, consider a square FGM plate made up of Aluminum-Alumina with $a/h=100$. \fref{fig:flutFGMRes} shows the influence of the boundary conditions on the critical aerodynamic pressure for various gradient index. It can be seen that the critical pressure is more for the clamped plate in comparison with that of the simply supported plate as expected. It is also seen that the aerodynamic pressure decreases with increase in the gradient index $n$. However the rate of decrease is high for low values of $n$. This can be attributed to the fact that the stiffness is high for the ceramic plate and minimum for the metallic plate and it degrades gradually with increase in the gradient index $n$. 

\section{Conclusions}
In this paper, we applied the NURBS based Bubnov-Galerkin iso-geometric finite element method  to study the static and dynamic response of functionally graded material plates. The first order shear deformation plate theory (FSDT) was used to describe the plate kinematics. Of course the present method is not limited to  FSDT and can easily be extended to higher order plate theories. It is to be noted that with NURBS basis functions, geometry could be exactly represented. Although in the present study only simple geometries are considered, the only thing that would change is the information pertaining to the geometry represented by the NURBS basis functions, when it is applied to model and/or analyze complex geometries.  The formulation when applied to thin plates, suffers from shear locking, which is alleviated by employing a modified shear correction factor. Numerical experiments have been conducted to bring out the influence of the gradient index, the plate aspect ratio and the plate thickness on the global response of functionally graded material plates. From the detailed numerical study, it can be concluded that with increasing gradient index $n$, the static deflection increases, while the free flexural vibration, critical buckling load and the flutter frequency decreases. This can be attributed to the reduction in stiffness of the material structure due to increase in the metallic volume fraction.

{\bf Acknowledgements} \\
The financial support of European Marie Curie Initial Training Network (FP7-People programme) is gratefully acknowledged.

\bibliographystyle{elsarticle-num}
\bibliography{nurbsPlate}

\begin{thebibliography}{10}
\expandafter\ifx\csname url\endcsname\relax
  \def\url#1{\texttt{#1}}\fi
\expandafter\ifx\csname urlprefix\endcsname\relax\def\urlprefix{URL }\fi
\expandafter\ifx\csname href\endcsname\relax
  \def\href#1#2{#2} \def\path#1{#1}\fi

\bibitem{koizumi1993}
M.~Koizumi, The concept of {FGM}, Ceramic Transactions - Functionally graded
  materials 34 (1993) 3--10.

\bibitem{Reddy2000}
J.~N. Reddy, Analysis of functionally graded plates, International Journal for
  Numerical Methods in Engineering 47 (2000) 663--684.

\bibitem{Yang2002}
J.~Yang, H.~S. Shen, Vibration characteristic and transient response of
  shear-deformable functionally graded plates in thermal environments, Journal
  of Sound and Vibration 255 (2002) 579--602.

\bibitem{Sundararajan2005}
N.~Sundararajan, T.~Prakash, M.~Ganapathi, Nonlinear free flexural vibrations
  of functionally graded rectangular and skew plates under thermal
  environments, Finite Elements in Analysis and Design 42~(2) (2005) 152--168.

\bibitem{Qian2004a}
L.~C. Qian, R.~C. Batra, L.~M. Chen, Static and dynamic deformations of thick
  functionally graded elastic plates by using higher order shear and normal
  deformable plate theory and meshless local {P}etrov {G}alerkin method,
  Composites Part B: Engineering 35 (2004) 685--697.

\bibitem{Ferreira2006}
A.~J.~M. Ferreira, R.~C. Batra, C.~M.~C. Roque, L.~K. Qian, R.~M.~N. Jorge,
  Natural frequencies of functionally graded plates by a meshless method,
  Composite Structures 75 (2006) 593--600.

\bibitem{natarajanmanickam2012}
S.~Natarajan, G.~Manickam, Bending and vibration of functionally graded
  material sandwich plates using an accurate theory, Finite Elements in
  Analysis and Design 57 (2012) 32--42.

\bibitem{ganapathiprakash2006}
M.~Ganapathi, T.~Prakash, N.~Sundararajan, Influence of functionally graded
  material on buckling of skew plates under mechanical loads, {ASCE} Journal of
  Engineering Mechanics 132 (2006) 902--905.

\bibitem{bathedvorkin1985}
K.~J. Bathe, E.~Dvorkin, A four node plate bending element based on {M}indlin -
  {R}eissner plate theory and mixed interpolation., International Journal for
  Numerical Methods in Engineering 21 (1985) 367--383.

\bibitem{somashekarprathap1987}
B.~R. Somashekar, G.~Prathap, C.~R. Babu, A field-consistent four-noded
  laminated anisotropic plate/shell element, Computers and Structures 25 (1987)
  345--353.

\bibitem{ganapathivaradan1991}
M.~Ganapathi, T.~K. Varadan, B.~S. Sarma, Nonlinear flexural vibrations of
  laminated orthotropic plates, Computers and Structures 39 (1991) 685--688.

\bibitem{bletzingerbischoff2000}
K.~U. Bletzinger, M.~Bischoff, E.~Ramm, A unified approach for shear-locking
  free triangular and rectangular shell finite elements, International Journal
  for Numerical Methods in Engineering 75 (2000) 321--334.

\bibitem{wangchen2004}
D.~Wang, J.~S. Chen, Locking-free stabilized conforming nodal integration for
  mesh-free {M}indlin-{R}eissner plate formulation, Computer Methods in Applied
  Mechanical and Engineering 193 (2004) 1065--1083.

\bibitem{nguyenrabczuk2008}
N.~T. Nguyen, T.~Rabczuk, H.~Nguyen-Xuan, S.~Bordas, A smoothed finite element
  method for shell analysis, Computer Methods in Applied Mechanics and
  Engineering 198 (2008) 165--177.

\bibitem{nguyen-xuanrabczuk2008}
H.~Nguyen-Xuan, T.~Rabczuk, S.~Bordas, J.~F. Debongnie, A smoothed finite
  element method for plate analysis, Computer Methods in Applied Mechanics and
  Engineering 197 (2008) 1184--1203.

\bibitem{kanok-nukulchaibarry2001}
W.~Kanok-Nukulchai, W.~Barry, K.~Saran-Yasoontorn, P.~H. Bouillard, On
  elimination of shear locking in the element-free {G}alerkin method,
  International Journal for Numerical Methods in Engineering 52 (2001)
  705--725.

\bibitem{He2001}
X.~Q. He, T.~Y. Ng, S.~Sivashanker, K.~M. Liew, Active control of {FGM} plates
  with integrated piezoelectric sensors and actuators, International Journal of
  Solids and Structures 38 (2001) 1641--1655.

\bibitem{Liew1994}
K.~M. Liew, K.~C. Hung, K.~M. Lim, A solution method for analysis of cracked
  plates under vibration., Engineering fracture mechanics 48~(3) (1994)
  393--404.

\bibitem{Ng2000}
T.~Y. Ng, K.~Y. Lam, K.~M. Liew, Effect of {FGM} materials on parametric
  response of plate structures, Computer Methods in Applied Mechanics and
  Engineering 190 (2000) 953--962.

\bibitem{Yang2001}
J.~Yang, H.~S. Shen, Dynamic response of initially stressed functionally graded
  rectangular thin plates, Composite Structures 54 (2001) 497--508.

\bibitem{matsunaga2008}
H.~Matsunaga, Free vibration and stability of functionally graded plates
  according to a 2{D} higher-order deformation theory, Composite Structures 82
  (2008) 499--512.

\bibitem{Vel2002}
S.~S. Vel, R.~C. Batra, Exact solutions for thermoelastic deformations of
  functionally graded thick rectangular plates, AIAA J 40 (2002) 1421--1433.

\bibitem{Vel2004}
S.~S. Vel, R.~C. Batra, Three-dimensional exact solution for the vibration of
  functionally graded rectangular plates, Journal of Sound and Vibration 272
  (2004) 703--730.

\bibitem{birman1995}
V.~Birman, Buckling of functionally graded hybrid composite plates, in:
  Proceedings of 10$^{th}$ Conference on Engineering Mechanics, Vol.~2, 1995,
  pp. 1199--1202.

\bibitem{javaherieslami2002}
R.~Javaheri, M.~Eslami, Buckling of functionally graded plates under in-plane
  compressive loading, ZAMM 82 (2002) 277--283.

\bibitem{woomeguid2003}
J.~Woo, S.~Meguid, K.~Liew, Thermomechanical postbuckling analysis of
  functionally graded plates and shallow cylindrical shells, Acta Mech. 165
  (2003) 99--115.

\bibitem{Prakash2006}
T.~Prakash, M.~Ganapathi, Supersonic flutter characteristics of functionally
  graded flat panels including thermal effects, Composite Structures 72 (2006)
  10--18.

\bibitem{Haddadpour2007}
H.~Haddadpour, H.~Navazi, F.~Shadmehri, Nonlinear oscillations of a fluttering
  functionally graded plate, Composite Structures 79 (2007) 242--250.

\bibitem{Sohn2008}
K.-J. Sohn, J.-H. Kim, Structural stability of functionally graded panels
  subjected to aero-thermal loads, Composite Structures 82 (2008) 317--325.

\bibitem{Sohn2009}
K.-J. Sohn, J.~Kim, Nonlinear thermal flutter of functionally graded panels
  under a supersonic flow, Composite Structures 88 (2009) 380--387.

\bibitem{Yang2010}
J.~Yang, Y.~Hao, W.~Zhang, S.~Kitipornchai, Nonlinear dynamic response of a
  functionally graded plate with a through-width surface crack, Nonlinear
  {D}ynamics 59 (2010) 207--219.

\bibitem{Kitipornchai2009}
S.~Kitipornchai, L.~Ke, J.~Y. andY Xiang, Nonlinear vibration of edge cracked
  functionally graded {T}imoshenko beams, Journal of Sound and Vibration 324
  (2009) 962--982.

\bibitem{natarajanbaiz2011}
S.~Natarajan, P.~Baiz, M.~Ganapathi, P.~Kerfriden, S.~Bordas, Linear free
  flexural vibration of cracked functionally graded plates in thermal
  environment, Computers and Structures 89 (2011) 1535--1546.

\bibitem{natarajanbaiz2011a}
S.~Natarajan, P.~Baiz, S.~Bordas, P.~Kerfriden, T.~Rabczuk, Natural frequencies
  of cracked functionally graded material plates by the extended finite element
  method, Composite Structures 93 (2011) 3082--3092.

\bibitem{baiznatarajan2011}
P.~M. Baiz, S.~Natarajan, S.~Bordas, P.~Kerfriden, T.~Rabczuk, Linear buckling
  analysis of cracked plates by {SFEM and XFEM}, Journal of Mechanics of
  Materials and Structure 6 (2011) 1213--1238.

\bibitem{veigabuffa2012}
L.~B.~a. da~Veiga, A.~Buffa, C.~Lovadina, M.~Martinelli, G.~Sangalli, An
  iso-geometric method for the reissner-mindlin plate bending problem, Computer
  Methods in Applied Mechanics and Engineering 209--212 (2012) 45--53.

\bibitem{nukulchaibarry2001}
W.~Kanok-Nukulchai, W.~Barry, K.~Saran-Yasontorn, P.~Bouillard, On elimination
  of shear locking in the element-free galerkin method, International Journal
  for Numerical Methods in Engineering 52 (2001) 705--725.

\bibitem{kikuchiishii1999}
F.~Kikuchi, K.~Ishii, An improved 4-node quadrilateral plate bending element of
  the reissne-mindlin type, Computational Mechanics 23 (1999) 240--249.

\bibitem{wu2004}
L.~Wu, Thermal buckling of a simply supported moderately thick rectangular
  {FGM} plate, Composite Structures 64 (2004) 211--218.

\bibitem{Rajasekaran1973}
S.~Rajasekaran, D.~Murray, Incremental finite element matrices, {ASCE} Journal
  of Structural Divison 99 (1973) 2423--2438.

\bibitem{Birman1990}
V.~Birman, L.~Librescu, Supersonic flutter of shear deformation laminated flat
  panel, Journal of Sound and Vibration 139 (1990) 265--275.

\bibitem{Ganapathi1996}
M.~Ganapathi, M.~Touratier, Supersonic flutter analysis of thermally stressed
  laminated composite flat panels, Composite Structures 34 (1996) 241--248.

\bibitem{cottrellhughes}
J.~A. Cottrell, T.~J. Hughes, Y.~Bazilevs, Isogeometric analysis: {T}oward
  integration of {CAD} and {FEA}, John Wiley, 2009.

\bibitem{vinhphusimpson2012}
V.~P. Nguyen, R.~N. Simpson, S.~P. Bordas, T.~Rabczuk, An introduction to
  {I}sogeometric analysis with {MATLAB} implementation: {FEM and XFEM}
  formulations, in review.

\bibitem{singhaprakash2011}
M.~Singha, T.~Prakash, M.~Ganapathi, Finite element analysis of functionally
  graded plates under transverse load, Finite elements in Analysis and Design
  47 (2011) 453--460.

\bibitem{gilhooleybatra2007}
D.~Gilhooley, R.~Batra, J.~Xiao, M.~McCarthy, J.~Gillespie, Analysis of thick
  functionally graded plates by using higher order shear and normal deformable
  plate theory and {MLPG} method with radial basis functions, Composite
  Structures 80 (2007) 539--552.

\bibitem{leezhao2009}
Y.~Lee, X.~Zhao, K.~Liew, Thermo-elastic analysis of functionally graded plates
  using the element free $kp-${R}itz method, Smart Materials and Structures 18
  (2009) 035007.

\bibitem{nguyen-xuantran2012}
H.~Nguyen-Xuan, L.~V. Tran, H.~Thai, T.~Nguyen-Thoi, Analysis of functionally
  graded plates by an efficient finite element method with node-based strain
  smoothing, Thin Walled Structures 54 (2012) 1--18.

\bibitem{thainguyen-xuan2012}
C.~Thai, H.~Nguyen-Xuan, N.~Nguyen-Thanh, T.-H. Le, T.~Nguyen-Thoi, T.~Rabczuk,
  Static, free vibration and buckling analysis of laminated composite
  {Reissner-Mindlin} plates using {NURBS} based isogeometric approach,
  International Journal for Numerical Methods in Engineering 91 (2012)
  571--603.

\bibitem{crocevenini2007}
L.~Croce, P.~Venini, Finite elements for functionally graded {Reissner-Mindlin}
  plates, Computer Methods in Applied Mechanics and Engineering 193 (2007)
  705--725.

\bibitem{zhaoliew2009}
X.~Zhao, K.~Liew, Geometrically nonlinear analysis of functionally graded
  plates using the element free $kp-${Ritz} method, Computer Methods in Applied
  Mechanics and Engineering 198 (2009) 2796--2811.

\bibitem{nguyen-xuantran2011}
H.~Nguyen-Xuan, L.~V. Tran, T.~Nguyen-Thoi, H.~Vu-Do, Analysis of functionally
  graded plates using an edge-based smoothed finite element method, Composite
  Structures 93 (2011) 3019--3039.

\bibitem{zhaolee2009}
X.~Zhao, Y.~Lee, K.~Liew, Free vibration analysis of functionally graded plates
  using the element free $kp-${Ritz} method, Journal of Sound and Vibration 319
  (2009) 918--939.

\bibitem{hashemifadaee2011}
S.~H. Hashemi, M.~Fadaee, S.~Atashipour, A new exact analytical approach for
  free vibration of reissne-mindlin functionally graded rectangular plates,
  International Journal of Mechanical Sciences 53 (2011) 11--22.

\bibitem{ganapathiprakash2006a}
M.~Ganapathi, T.~Prakash, Thermal buckling of simply supported functionally
  graded skew plates, Composite Structures 74 (2006) 247--250.

\end{thebibliography}

\end{document}